\documentclass[11pt]{amsart}
\usepackage{vmargin}
\usepackage{psfig}
\setpapersize{A4}
\setmarginsrb{25mm}{25mm}{25mm}{25mm}{12pt}{10mm}{12pt}{10mm}

\newtheorem{theorem}{Theorem}[section]

\newtheorem{lemma}[theorem]{Lemma}

\newtheorem{proposition}[theorem]{Proposition}

\newtheorem{corollary}[theorem]{Corollary}

\theoremstyle{definition} 
\newtheorem{definition}[theorem]{Definition}
\newtheorem{example}[theorem]{Example}

\theoremstyle{remark} 

\newtheorem*{unremark}{Remark}

\makeatletter 
\def\romenumi{
\def\theenumi{\roman{enumi}}
\def\p@enumi{\theenumi}
\def\labelenumi{(\@roman\c@enumi)}}
\makeatother

\newif\ifShowLabels 
 \ShowLabelsfalse
\newdimen\mgheight \def\marginnotes#1{%
\leavevmode\vadjust{\setbox0=\hbox{{\tt \quad\quad {\small\rm #1}}}%
\mgheight=\ht0 \advance\mgheight by \dp0 \advance\mgheight by \lineskip
\kern -\mgheight \vbox to \mgheight{\rightline{\rlap{\box0}} \vss}}}
\def\tlabel#1{\label{#1} \ifShowLabels \marginnotes{#1} \fi}

\newcommand{\noproof}{\hfill\qedsymbol}

\DeclareMathOperator{\re}{Re}

\DeclareMathOperator{\Hess}{{\b{H}}}
 
\DeclareMathOperator{\rk}{rank}
\DeclareMathOperator{\cone}{\b{K}}
\DeclareMathOperator{\dir}{\b{dir}}
\DeclareMathOperator{\grad}{\nabla}

\def\sing{\mathcal{V}} 
\def\torus{\b{T}} 
\def\disk{\b{D}}
\def\expsmall{{\rm exponentially\ small}}
\def\poly{{\rm polynomial}}
\def\cmat{{\b{C}}}
\def\vv{{\bf v}}
\def\xx{{\bf x}}
\def\rr{{\bf r}}
\def\numer{{G}}
\def\denom{{H}} 
\def\mult{{\psi}} 
\def\simp{{\Delta}} 
\def\meas{{\mu}}
\def\nice{\mathcal{E}}
\def\rank{{\rm rank}}

\def\b#1{\mathbf{#1}} 
\def\bb#1{\boldsymbol{#1}} 
\def\w#1{\widehat{#1}}
 
\def\zh{\w{\b{z}}} 
\def\zsh{\w{\b{z^*}}}

\def\nbd{{\mathcal N}} 
\def\direc{{\bb{\delta}}} 
\def\domain{{\mathcal{D}}} 
\def\logdom{{\rm log} \domain} 
\def\rate{{\gamma}}
\def\hyper{\mathcal{H}} 
\def\weier{\mbox{multiple point }}
\def\weierdef{\mbox{\it multiple point}} 
 
\def\weierformdef{\mbox{\it factored form}} 
\def\R{{\mathbb{R}}}
\def\CC{{\mathbb{C}}}

\def\statset{{\mathcal{S}}}
\def\interior{\mbox{boundedly interior }}

\swapnumbers 
\numberwithin{equation}{section}

\begin{document} 
\title[Asymptotics of multivariate sequences II] {Asymptotics of
multivariate sequences {II}.\\Multiple points of the singular variety.}

\author{Robin~Pemantle} 
\address{Department of Mathematics, Ohio State University, Columbus OH
43210, USA}
\email{pemantle@math.upenn.edu}
\curraddr{Department of Mathematics, University of Pennsylvania,
Philadelphia, PA 19104-6395, USA}  
\thanks{Research supported in part by NSF grants DMS-9996406 and
DMS-0103635}

\author{Mark~C.~Wilson} 
\address{Department of Computer Science, University of Auckland,
Private Bag 92019 Auckland, New Zealand} 
\email{mcw@cs.auckland.ac.nz}

\date{\today} 

\begin{abstract}Let $F(\b{z})=\sum_\b{r} a_\b{r}\b{z^r}$ be a
multivariate  generating function which is meromorphic in some
neighborhood of the origin of $\mathbb{C}^d$, and let $\sing$ be its
set of singularities. Effective asymptotic expansions for the
coefficients can be obtained by complex contour integration near
points of $\sing$.

In the first article in this series, we treated the case of smooth points  of $\sing$. In this article we deal with multiple points of $\sing$. 
Our results show that the central limit (Ornstein-Zernike) behavior 
typical of the smooth case does not hold in the multiple point case.
For example, when $\sing$ has a multiple point singularity at 
$(1 , \ldots , 1)$, rather than $a_\b{r}$ decaying as 
$|\b{r}|^{-1/2}$ as $|\b{r}| \to \infty$, $a_\b{r}$ is very nearly
polynomial in a cone of directions.

\end{abstract}

\maketitle

\section{Introduction} \tlabel{sec:intro} 
In \cite{pemantle-wilson-multivariate1}, we began a series of articles
addressing the general problem  of computing asymptotic expansions for a
multivariate sequence whose generating function is known. Such problems
are encountered frequently in combinatorics and probability; see for
instance Examples~2~--~8 in Section~1 of~\cite{pemantle-lecnotes}, which
collects examples from various  sources
including~\cite{larsen-lyons-coalescing,flatto-mckean-queues-parallel,
wilf-GFology,comtet-advanced}. Our aim is to present methods which are as general as possible, and lead to effective computation.  Our apparatus may be applied to any function whose dominant singularities are poles. Among other things, we showed in~\cite{pemantle-wilson-multivariate1} that for all nonnegative bivariate sequences of this type, our method is  applicable (further work may be required in one degenerate case). 

The article \cite{pemantle-wilson-multivariate1} handled the case when
the pole variety was smooth.  The present work is motivated by
a collection of applications where the pole variety is composed of
intersecting branches.  Several examples are as follows.  
\begin{enumerate}
\item $$F(x,y) = \frac{2}{(1-2x)(1-2y)} \left[ 2+ \frac{xy-1}{1-xy(1+x+y+2xy)}\right]$$
arises in Markov modeling \cite{karloff};
\item $$F(\b{z}) = \frac{G(\b{z})}{\prod_{j=1}^k L_j (\b{z})}$$
where the $L_j$ are affine, arises in queueing theory \cite{bertozzi-mckenna-queueing};
\item $$F(x,y,z) = \frac {x^a y^b z^c}{(1-x)(1-y)(1-z)(1-xz)(1-yz)}$$
is a typical generating function arising in enumeration of integer
solutions to unimodular linear equations, and is the running example
used in \cite{deloera-sturmfels};
\item $$F(x,y) = \frac{1}{(1 - (1/3) x - (2/3) y)(1 - (2/3) x - (1/3) y)}$$
counts winning plays in a dice game (see Example~\ref{eg:comb}).
\item $$F(x,y,z) = \frac{z/2}{(1-yz) P(x,y,z)}$$
arises in the analysis of random lattice tilings \cite{cohn-elkies-propp-aztec}.  In this case, in addition to the factorization of the denominator, there is an isolated singularity in the variety defined by the polynomial $P$, for
which reason the analysis is not carried out in the present paper 
but in \cite{cohn-pemantle-fixation}.
\end{enumerate}
These examples serve as our motivation to undertake a categorical
examination of generating functions with self-intersections of
the pole variety.  As in the above cases, one may imagine that this
comes about from a factorization of the generating function, although
we show in Example~\ref{eg:figure 8} that the method works as well for 
irreducible functions with the same ``multiple point'' local geometry.

Our approach is analytic. For simplicity, we restrict to the two-variable
case in this introduction, though our methods work for any number of
variables. Given a sequence $a_{rs}$ indexed by the $2$-dimensional 
nonnegative integer lattice, we seek asymptotics as $r,s \to \infty$. 
Form the  generating function $F(z,w)=\sum_{r,s}a_{rs}z^rw^s$ of the
sequence; we assume that $F$ is analytic in some neighborhood of the
origin.

The iterated Cauchy integral formula yields
$$a_{rs}=\frac{1}{(2\pi i)^2} \int_{\mathcal{C}'} \int_{\mathcal{C}}
\frac{F(z,w)}{z^{r+1}w^{s+1}}\,dw \,dz,$$
where $\mathcal{C}$ and $\mathcal{C'}$ are circles centered at $0$ and
$F$ is analytic on a polydisk containing the torus $\mathcal{C}\times
\mathcal{C'}$. Expand the torus, by expanding (say) $\mathcal{C}$,
slightly beyond a {\em minimal singularity} of $F$ (that is, a point
$(z_0 , w_0)$ at which the expanding torus first touches the singular set
$\sing$ of $F$). The difference between the corresponding inner integrals
is then computed via residue theory.  Thus $a_{rs}$ is represented as a
sum of a residue and an integral on a large torus.  One hopes that the
residue term is dominant, and gives a good approximation to $a_{rs}$. 
The residue term is itself an integral (since the residue is taken in the
inner integral only).  A stationary phase analysis of the residue
integral shows when our hope is realized, and yields an asymptotic
expansion for $a_{rs}$.

In \cite{pemantle-wilson-multivariate1} we considered the case when the minimal
singularity in question is a smooth point of $\sing$.  The present article deals
with the case where the minimal singularity is a {\em multiple
point}: locally, the singular set is a union of finitely many graphs of analytic
functions. This case includes the smooth point analysis of \cite{pemantle-wilson-multivariate1}, and Theorem~\ref{thm:transverse} with $n$ set to $0$ essentially generalizes the analysis of that article). In the two-variable case, we know from Lemma~6.1 of \cite{pemantle-wilson-multivariate1}
that every minimal singularity must have this form or else be a cusp (see
also~\cite[Lemma~3.1] {tsikh-meromorphic-2D}), though more complicated
singularities may arise in higher dimensions.

The rest of this article is organized as follows.  In the remainder of this
section we describe in more detail the program begun
in~\cite{pemantle-wilson-multivariate1} and continued here and in future articles
in this series. Section~\ref{sec:prelim} deals with notation and preliminaries
required for the statement of our main results. Those results, along with
illustrative examples, are listed in Section~\ref{sec:results}.  Proofs are given
in Section~\ref{sec:proofs}.  We discuss some further details  and outline future
work in Section~\ref{sec:future}.

\subsection*{Details of the Program}

Our notation is similar to that in \cite{pemantle-wilson-multivariate1}.
For clarity we shall reserve the names of several objects throughout.  We
use {\bf boldface} to indicate a (row or column) vector.  The number of
variables will be denoted $d+1$.  The usual multi-index notation is in
use: $\b{z}$ denotes a vector $(z_1, \dots ,z_{d+1})^T \in
\mathbb{C}^{d+1}$, and if $\b{r}$ is an integer vector then
$\b{z}^{\b{r}}=\prod_j z_j^{r_j}$. We also use the convention that a
function, ostensibly of $1$ variable, applied to an element of 
$\mathbb{C}^{d+1}$ acts on each coordinate separately --- for example
$e^{\b{x}}=(e^{x_1}, \dots , e^{x_{d+1}}).$  Throughout, $\numer$ and
$\denom$ denote functions analytic in some polydisk about $\b{0}$ and
$F=G/H=\sum_{\b{r}} a_{\b{r}}\b{z}^\b{r}$. The set where $H$ vanishes
will be called the {\em singular variety} of $F$ and denoted by $\sing$. 
Let $\disk (\b{z}), \torus (\b{z})$ denote respectively the polydisk and
torus (both centered at the origin) on which $\b{z}$ lies.

A crude preliminary step in approximating $a_\b{r}$ is to determine its
exponential rate; in other words, to estimate $\log |a_\b{r}|$ up to a
factor of $1 + o(1)$.  Let $\domain$ denote the (open) domain of
convergence of $F$ and let $\logdom$ denote the logarithmic domain in
$\R^{d+1}$, that is, the set of $\b{x} \in \R^{d+1}$ such that
$e^{\b{x}} \in \domain$.  If $\b{z^*} \in \domain$ then Cauchy's integral
formula 

\begin{equation} \tlabel{eq:cauchy-multi} a_\b{r} = \left (
\frac{1}{2 \pi i} \right )^{d+1} \int_{\torus (\b{z^*})}
\frac{F(\b{z})}{\b{z}^{\b{r}+\b{1}}} \, d\b{z} 
\end{equation} 
shows that $a_\b{r} = O(|\b{z^*}|^{-\b{r}})$.  Letting $\b{z^*} \to \partial
\domain$ gives $$\log | a_\b{r} | \leq - \b{r} \cdot \log |\b{z^*}| +
o(|\b{r}|),$$ and optimizing in $\b{z^*}$ gives $\log |a_\b{r} | \leq \rate
(\b{r}) + o(|\b{r}|)$ where
 
\begin{equation} \tlabel{eq:rate}
\rate (\b{r}) := - \sup_{\b{x} \in \logdom} \b{r} \cdot \b{x} \, .
\end{equation}

The cases in which the most is known about $a_\b{r}$ are those in which
this upper bound is correct, that is, $\log |a_\b{r}| = \rate (\b{r}) +
o(|\b{r}|)$.  To explain this, note first that the supremum
in~(\ref{eq:rate}) is equal to $\b{r} \cdot \b{x}$ for some $\b{x} \in
\partial \logdom$.  The torus $\torus (\b{e^x})$ must contain some minimal
singularity $\b{z^*} \in \sing \cap \partial \domain$.  Asking that $\log
|a_\b{r}| \sim - \b{r} \cdot \b{z^*}$ is then precisely the same as
requiring the Cauchy integral~(\ref{eq:cauchy-multi}) --- or the residue
integral mentioned above --- to be of roughly the same order
as its integrand.  This is the situation in which it easiest to estimate
the integral.

Our program may now be summarized as follows.  Associated to each minimal singularity $\b{z^*}$ is a cone $\cone (\b{z^*}) \subseteq (\R^+)^{d+1}$. 
Given $\b{r}$, we find one or more $\b{z^*} = \b{z^*} (\b{r}) \in \sing \cap
\partial \domain$ where the upper bound is least.  We then attempt to
compute a residue integral there.  This works only if $\b{r} \in \cone
(\b{z^*})$ and if the residue computation is of a type we can handle.  Our
program is guaranteed to succeed in some cases, and conjectured to
succeed in others.  It is known to fail only in some cases where the
$a_\b{r}$ are not nonnegative reals (not the most important case in
combinatorial or probabilistic applications) and even then a variant
seems to work.

To amplify on this, define a point $\b{z^*} \in \sing$ to be {\em minimal}
if $\b{z^*} \in \partial \domain$ and each coordinate of $\b{z^*}$ is nonzero. There are only three possible types of minimal singularities~\cite[Lemma~6.1]
{pemantle-wilson-multivariate1}), namely smooth points of $\sing$,
multiple points and cone points (all defined below). It is conjectured
that for all three types of points, and any $\b{r} \in \cone (\b{z^*})$, we indeed have 

$$
\log |a_\b{r}| = \rate (\b{r}) + o(|\b{r}|) = - \b{r} \cdot \log
|\b{z^*}| + o(|\b{r}|) \, .  
$$ 
This is proved for smooth points
in~\cite{pemantle-wilson-multivariate1} via residue integration, and the
complete asymptotic series obtained.  It is proved in the present work
for multiple points under various assumptions; the fact that
these do not cover all cases seems due more to taxonomical problems
rather than the inapplicability of the method. The problem remains open
for cone points, along with the problem of computing asymptotics.

It is also shown  in~\cite{pemantle-wilson-multivariate1} that when $a_\b{r}$ are all nonnegative, then $\b{r} \in \cone (\b{z^*} (\b{r}))$, and therefore that a resolution of the above problem yields a complete
analysis of nonnegative sequences with pole singularities.

When the hypothesis that $a_\b{r}$ be real and nonnegative is removed, it is not always true that $\b{r}\in \cone (\b{z^*} (\b{r}))$.  In all
examples we have worked, there have been points $\b{z^*} \notin
\overline{\domain}$  for which $\b{r} \in \cone(\b{z^*})$ and a residue
integral near $\b{z^*}$ may be proved a good approximation to $a_\b{r}$;
the method is to contract the torus of integration to the origin in some
way other than simply iterating the contraction of each coordinate circle to a point. Thus a second open question is to settle whether there is always such a point $\b{z^*}$ in the case of mixed signs (see the
discussion in~\cite[Section~7]{pemantle-wilson-multivariate1}).

The scope of the present article is as follows. We define multiple
points and carry out the residue integral arising near a
strictly minimal, multiple point.  Unlike the case for smooth points, the
integral is not readily recognizable as a standard multivariate
Fourier-Laplace integral, and a key result of this article is a more
manageable representation of the residue to be integrated
(Corollary~\ref{cor:residue sum} and Lemma~\ref{lem:tooscint}). 
The Fourier-Laplace integrals arising fall just outside the scope of the standard references, and we have been led to develop generalizations of known results. These results, which would take too much space here, are included in \cite{pemantle-wilson-analysis}.

The asymptotics arising from multiple point singularities differ
substantially from asymptotics in the smooth case.  In the remainder of
this introduction, we give examples to illustrate this.  

\begin{figure}
\centerline{\psfig{figure=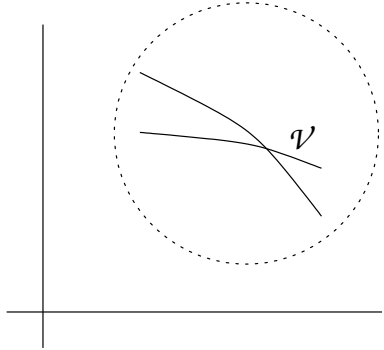,width=2in}}
\vspace{10pt}
\caption{Local picture of a double pole}
\tlabel{fig:simple}
\end{figure}

\begin{example}[simplest possible multiple point] \tlabel{eg:simple} 
Let $F(z,w)$ be a two-variable generating function and suppose that  the point
$(1,1)$ is a double pole of $F$ (thus $F = \numer /  \denom$ with $\numer (z,w)
\neq 0$ and $\denom (z,w)$ vanishing to order $2$ at $(1,1)$).  If $F$ has no other poles $(z,w)$ with $|z| , |w| \leq 1$, and if the two branches of the singular variety $\sing$ meet transversely at $(1,1)$ as in Figure~\ref{fig:simple}, then for some positive constants $c$ and $C$,

\begin{equation} \tlabel{eq:plateau 1}
a_{rs} = C + O(e^{-c |(r,s)|})
\end{equation}
for all $(r,s)$ in a certain cone, $\cone$, in the positive integer quadrant.  
This is proved in Theorem~\ref{thm:double-point-generic} below, 
and the constant $C$ computed.  Exact statement of the transversality
hypothesis requires some discussion of the geometry of $\sing$.  The
constant $C$ is computed in terms of some algebraic quantities
derived from $F$.  The cone $\cone$ in which this holds is easily described
in terms of the tangents to the branches of the double pole.
The need for some preliminary algebraic and geometric analysis 
to define transversality and to compute $C$ and $\cone$
motivates our somewhat lengthy Section~\ref{sec:prelim}. 

Example~\ref{eg:figure 8} below shows the details of this computation 
for a particular $F$.  The exponential bound on the error follows
from results in~\cite{pemantle-finite-dimension} which are cited in
Section~\ref{sec:proofs}.  If the multiple pole is moved to a point $(z^*,w^*)$ 
other than $(1,1)$, a factor of $(z^*)^{-r} (w^*)^{-s}$ is introduced.  

Compare this with the case where $\sing$ is smooth, intersecting the
positive real quadrant as shown in Figure~\ref{fig:smooth}.
In this case one has Ornstein-Zernike (central limit) behavior. Suppose, as above, that $(1,1)\in \sing$ and  $F$ has no other poles
$(z,w)$ with $|z| , |w| \leq 1$. Then 
as $(r,s) \to \infty$ with $r/s$ fixed, $a_{rs}$ is rapidly decreasing 
for all but one value of $r/s$; for that distinguished direction, 
$a_{rs} \sim C |(r,s)|^{-1/2}$ for some constant, $C$, and the error
terms may be developed in a series of decreasing powers of $|(r,s)|$.
The two main differences between the double and single pole cases
are thus the existence of a plateau in the double pole case, and the
flatness up to an exponentially small correction 
versus a correction of order $|(r,s)|^{-3/2}$ for a single pole.  
\noproof
\end{example}

\begin{figure}
\centerline{\psfig{figure=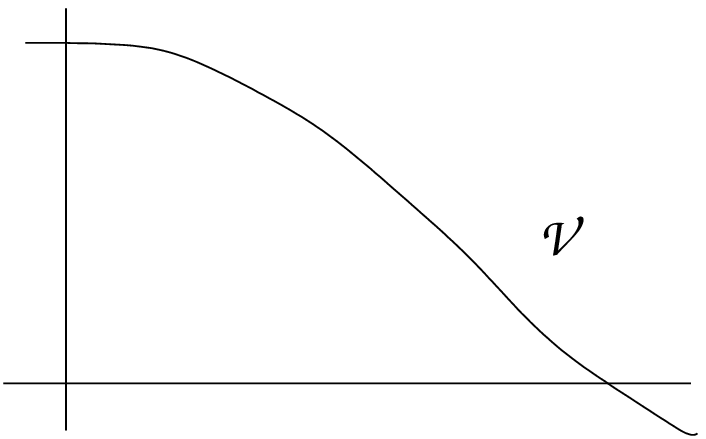,width=2in}}
\vspace{10pt}
\caption{When $\sing$ has only smooth points}
\tlabel{fig:smooth}
\end{figure}

\begin{example}[pole of greater order] \tlabel{ex:higher pole}
Alter the previous example so that $\sing$ has a pole of some 
order $n+1$ at $(1,1)$, as in Figure~\ref{fig:tent}.  
Then the formula~(\ref{eq:plateau 1}) becomes instead

$$
a_{rs} = P(r,s) + O(e^{-c |(r,s)|})
$$
where $P$ is a piecewise polynomial of degree $n-1$.  
In Example~\ref{eg:distinct-tangents}
below, $P$ is explicitly computed in the case $n+1 = 3$.
Again, we see the chief differences from the smooth case 
being a cone of non-exponential decay, and an exponentially small
correction to the polynomial $P(r,s)$ in the interior of the cone.
\noproof
\end{example}
\begin{figure}
\centerline{\psfig{figure=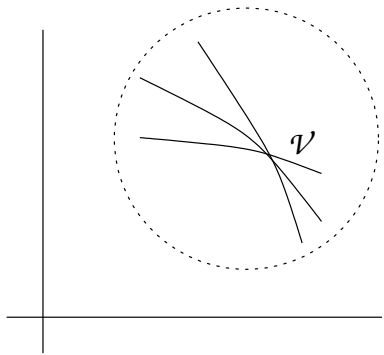,width=2in}}
\vspace{10pt}
\caption{A pole of order 3}
\tlabel{fig:tent}
\end{figure}

In higher dimensions there are more possible behaviors, but
the same sorts of results hold.  There is a cone on which the
exponential rate has a plateau; under some conditions the correction
terms within this cone are exponentially small.  The following table
summarizes the results proved in this paper; transversality
assumptions have been omitted. \hfill \\[2ex]

\begin{center}
\begin{tabular}{|l|l|l|}
\hline
{\bf Theorem number} & {\bf Hypotheses} & {\bf Asymptotic behavior} \\
\hline
\hline
Theorem~\ref{thm:double-point-generic} & two curves in 2-space &
   $C + \expsmall$ \\
\hline
Theorem~\ref{thm:compnondeg} & $d+1$ sheets in $(d+1)$-space &
   $C + \expsmall$ \\
\hline
Theorem~\ref{thm:nondeg} & more sheets than the dimension &
   $\poly + \expsmall$ \\
\hline
Theorem~\ref{thm:transverse} & fewer sheets than the dimension &
   asymptotics start with $|\rr|^{n/2-d/2}$ \\
\hline
Theorem~\ref{thm:degenerate} & all sheets tangent &
   asymptotics start with $|\rr|^{n-d/2}$ \\
\hline
\end{tabular} ~~\\[4ex]
\end{center}

We conclude the introduction by giving a combinatorial application of 
Example~\ref{eg:simple}.

\begin{example}[combinatorial application] \tlabel{eg:comb}

An independent sequence of random numbers uniform on $[0,1]$ is
used to generate biased coin-flips: if $p$ is the probability
of heads then a number $x \leq p$ means heads and $x > p$ means tails.
The coins will be biased so that $p=2/3$ for the first 
$n$ flips, and $p=1/3$ thereafter.
A player desires to get $r$ heads and $s$ tails and is allowed
to choose $n$.  On average, how many choices of $n \leq r + s$
will be winning choices?  

The probability that $n$ is a winning choice for the player 
is precisely
$$
\sum_{a+b=n}\binom{n} {a} (2/3)^a (1/3)^b \binom{r+s-n} {r-a} (1/3)^{r-a} 
(2/3)^{s-b} \, .
$$
Let $a_{rs}$ be this expression summed over $n$.  The array
$\{ a_{rs} \}_{r , s \geq 0}$ is just the convolution of the
arrays $\binom{r+s} {r} (2/3)^r (1/3)^s$ and $\binom{r+s} {r} (1/3)^r 
(2/3)^s$, so the generating function $F(z,w) := \sum a_{rs} z^r w^s$
is the product
$$
F (z,w) = \frac{1} {(1 - \frac{1} {3} z - \frac{2} {3} w) 
   (1 - \frac{2} {3} z - \frac{1} {3} w)} \, .
$$
Applying Theorem~\ref{thm:double-point-generic} 
with $G \equiv 1$ and $\det \Hess = -1/9$, 
we see that $a_{rs} = 3$ plus a correction which is exponentially
small as $r,s \to \infty$ with $r/(r+s)$ staying in any subinterval 
of $(1/3,2/3)$.  A purely combinatorial analysis of the sum 
may be carried out to yield the leading term, 3, but says nothing
about the correction terms.  The diagonal extraction method 
of~\cite{hautus-klarner-diagonal} yields very precise information for 
$r=s$ but nothing more general in the region $1/3< r/(r+s) < 2/3$.  
\noproof
\end{example}

\section{Preliminary definitions and notation} \tlabel{sec:prelim}

For each $\b{z}\in\mathbb{C}^{d+1}$, the truncation $(z_1, \dots ,z_d)$ will be
denoted $\w{\b{z}}$, and the last coordinate $z_{d+1}$ simply by $z$. We do not
specify the size for constant vectors --- for example $\b{1}$ denotes the vector
$(1, \dots ,1)^T$ of whatever size is appropriate. Thus we write
$\w{\b{1}}=\b{1}$.  A minimal singularity $\b{z^*}$ is {\em strictly minimal} if
$\sing \cap \disk (\b{z^*}) = \{ \b{z^*} \}$. When a minimal point is not strictly
minimal, one must add (or integrate) contributions from all points of $\sing \cap \torus (\b{z^*})$.  This step is routine and will be carried out in a future article; we streamline the exposition here by assuming strict minimality.

This article deals  entirely with minimal points $\b{z^*}$ of $\sing$ near
which $\sing$ decomposes as a union of sheets $\sing_j$, each of which is
a graph of an analytic function $z=u_j(\w{\b{z}})$. The algebraic
description of this situation is as follows.  By the Weierstrass
Preparation Theorem~\cite{griffiths-harris-principles} there is a
neighborhood of $\b{z^*}$ in which we may write $H(\b{z}) = \chi (\b{z})
W(\b{z})$ where $\chi$ is analytic and nonvanishing and $W$ is a {\em
Weierstrass polynomial}. This  means that $$W(\b{z}) = z^{n+1} +
\sum_{j=0}^n \chi_j (\zh) z^j $$ where the multiplicity $n+1$ is at least
2 if $\b{z^*}$ is not a smooth point, and the analytic functions $\chi_j$
vanish at $\zsh$.

Now suppose that $\b{z^*}$ has all coordinates nonzero. Recalling from
Lemma~6.1 of~\cite{pemantle-wilson-multivariate1} that any such minimal
point of $\sing$ is locally homogeneous, we see in fact that $\chi_j$
vanishes to homogeneous degree $n+1-j$ at $\zsh$ and that $\chi_0$ has
nonvanishing pure $z_i^{n+1}$ terms for each $1 \leq i \leq d$.  Then
$\sing$ is locally the union of smooth sheets if the degree $n+1$
homogeneous part of $W$ (the leading term) factors completely into 
{\em distinct} linear factors, while if $\sing$ is locally the union 
of smooth sheets then we have a factorization

\begin{equation}
\tlabel{eq:local fac} H(\b{z})= \chi (\b{z}) \prod_{j=0}^n [z - u_j
(\zh)] 
\end{equation} 
for analytic ({\em not necessarily distinct}) functions
$u_j$ mapping a neighborhood of $\zsh$ to a neighborhood of $z^*$.

\begin{unremark} If the expansion of $H$ near $\b{z^*}$ vanishes to order
$n+1$ in $z$, then there are always $n+1$ solutions (counting
multiplicity) to $H(\zh , z) = 0$ for $\zh$ near $\zsh$.  These vary
analytically, but may be parametrized by $n+1$ analytic functions $u_j$
only if there is no monodromy, that is, if one stays on the same branch
moving around a cycle in the complement of the singular set.  Thus a
\weier singularity is one whose monodromy group is trivial.  We remark
also that in two variables, local homogeneity of the minimal point
$\b{z^*}$ nearly implies it is a multiple point.  The only other possibility
is a cusp whose tangents are all equal; such a singularity turns out
not to affect the leading order asymptotics and to affect the lower
order asymptotics only in a single direction.
\end{unremark}

It turns out to be more convenient to deal with the reciprocals
$v_j=1/u_j$. The basic setup throughout the rest of this article is as
follows.

\begin{definition} \tlabel{defn:preparation} The point $\b{z^*}$ of
$\sing$ is a $\weierdef$ if there are analytic functions $v_0, \dots
,v_n, \phi$ and a local factorization 

\begin{equation}
\tlabel{eq:weier} F(\b{z})=\frac{\phi(\b{z})}{\prod_{j=0}^n
\left(1-zv_j(\w{\b{z}})\right)} \, ,
\end{equation} 
which we call the $\weierformdef$ of $F$, such that
 
\begin{itemize} \romenumi
\item 
$(z^*)_{d+1}v_j (\w{\b{z^*}})=1$ for all $j$; 
\item 
each $\displaystyle \frac{\partial v_j}{\partial z_k} (\w{\b{z^*}}) \neq 0$; 
\item
$z_{d+1}v_j (\w{\b{z}})=1$ for some $j$ if and only if $\b{z} \in \sing$. 
\end{itemize}
\end{definition}

\begin{unremark}

Let $\sing_j$ denote the local hypersurface parametrized by $z = u_j
(\zh)$. Then the first of the conditions says that each $\sing_j$ passes
through $\b{z^*}$: there are no extraneous factors in the denominator
representing surface elements not passing through $\b{z^*}$.  The last
condition says that the zeros of the denominator are exactly the poles of
$F$: there are no extraneous factors in the denominator vanishing at
$\b{z^*}$ and cancelling a similar divisor in the numerator.

It may appear that some generality has been lost in imposing the second
condition, since we are assuming that each sheet of $\sing$ projects
diffeomorphically onto any coordinate hyperplane. This latter property
is in fact guaranteed by the non-vanishing of the pure $z_i^{n+1}$
terms of $W$.  

To compute $\phi$ directly from $\numer$ and $\denom$, 
differentiate~(\ref{eq:local fac}) $n+1$ times in the $z := z_{d+1}$
coordinate at the point $\b{z^*}$ to write
$$
\left ( \frac{\partial} {\partial z} \right )^{n+1} 
\denom (\b{z^*}) = (n+1)! \, \chi (\b{z^*}) \, .
$$
We may then write
$$
\frac{\phi} {\prod_{j=0}^n (1 - z v_j \b{\w{z}})} = \frac{\numer} 
   {\denom} = \frac{\numer} {\chi \prod_{j=0}^n (-u_j) \prod_{j=0}^n
   (1 - z v_j (\b{\w{z}}))}
$$
and solve for $\phi$ at $\b{z} = \b{z^*}$ to obtain
\begin{equation} \tlabel{eq:comp phi}
\phi (\b{z^*}) = \frac{(n+1)!} {(-(z^*)_{d+1})^{n+1}} \, \frac{G(\b{z^*})} 
   {\left ( \frac{\partial} {\partial z_{d+1}} \right )^{n+1} 
   \denom (\b{z^*})} \, .
\end{equation}
\noproof
\end{unremark}

For the remaining definitions, fix a strictly minimal element
$\b{z^*}\in\mathcal{V}$ which is a multiple point, and let $v_0, \dots ,v_n$ be as above.

\begin{definition}[Cone of directions corresponding to multiple point]
\tlabel{defn:cone}
For each sheet $\sing_j$ and multiple point $\b{z^*}$, let $\dir_j(\b{z^*})$ be the  vector defined by  
$$\dir_j(\b{z^*}) := \left. \left(\frac{z_1}{z_{d+1}} \frac{\partial v_j}{\partial z_1}, \dots , \frac{z_d}{z_{d+1}}\frac{\partial v_j}{\partial z_d}, 1\right) \right|_{\b{z} = \b{z^*}}.$$
We denote by $\cone(\b{z^*})$ the positive hull of all the $\dir_j(\b{z^*})$ and  by $\cone_0(\b{z^*})$ their convex hull, in other words the intersection of $\cone(\b{z^*})$ with the hyperplane $z_{d+1}=1$. Geometrically, $\cone$ is precisely the collection of outward normal vectors to support hyperplanes of the logarithmic domain  of convergence of $F$ at the point $(\log |z^*_1| , \ldots , \log |z^*_{d+1}| )$; see~\cite{pemantle-wilson-multivariate1} for details.

Let $\cmat(\b{z^*})$ be the matrix whose $j$th row is 
$\dir_j(\b{z^*})$. We say that $\b{z^*}$ is {\bf
nondegenerate} if the rank of $\cmat$ is $d+1$, {\bf transverse} if the
rank is $n+1$, and {\bf completely nondegenerate} if it is both
transverse and nondegenerate; in this case necessarily $n=d$ and the
multiplicity of each sheet is $1$.
 
If $K$ is a subset of $\cone$ or $\cone_0$  consisting of
vectors whose directions are bounded away from the walls (equivalently
the image of $K$  under the natural map \medspace $\bar{}: \mathbb{R}^{d+1}
\setminus \{ \b{0} \} \to \mathbb{RP}^d$ is a compact subset  of the
interior of $\overline{\cone}$) then we shall say
that  $K$ is {\bf \interior} to $\cone$ (or $\cone_0$).
\noproof

\end{definition}

\begin{unremark}
The importance of $\cone$ is that analysis near $\b{z^*}$ will yield
asymptotics for $\b{r}$ in precisely the directions in $\cone(\b{z^*})$.
In the smooth case, $n = 0$ and so $\cone(\b{z^*})$ reduces to a single
ray. The point $\b{z^*}$ is transverse if and only if the normals there
to the surfaces $\sing_j$ span a space of dimension $n+1$.
Alternatively, the $n+1$ tangent hyperplanes intersect transversally.
When $n>d$, transversality must be violated; nondegeneracy means that
there is as little violation as possible.  
\noproof
\end{unremark}

\section{Main theorems and illustrative examples} \tlabel{sec:results}
All of the results in this section are ultimately proved by reducing
the problem to the computation of asymptotics for a Fourier-Laplace
integral (the proofs are presented in later sections). Owing to the
large number of possibilities arising in this analysis, constructing a
complete taxonomy of cases is rather challenging, and the number of 
potential theorems is enormous.

We have chosen to present a series of theorems of varying complexity and
generality. Taken together, they completely describe asymptotics
associated with nondegenerate multiple points. Our analysis of other
types of multiple points requires additional (mild) hypotheses, which
will almost always be satisfied in applications. The most important case
is that of transversal points, but we also treat various types of
tangencies and degeneracies. However, it is always possible that a
practical problem involving a meromorphic generating function may not fit
neatly into our classification scheme. We hope to convince the reader
that in such a situation our basic method will yield Fourier-Laplace
integrals from which asymptotics can almost certainly be extracted in a
systematic way.   

In two variables, our results specialize to the following cases. If 
$\sing$ is locally the union of $n+1$ analytic graphs, then either  
at least two tangents are distinct, in which case Theorem~\ref{thm:nondeg} applies,  or all tangents coincide. This latter case is more complicated and we require some extra hypotheses; see Theorem~\ref{thm:degenerate}.  

We begin with the simplest case, in which the result admits a
relatively self-contained statement, with as little extra notation as
possible. 

\begin{theorem}[2 curves meeting transversally in 2-space]
\tlabel{thm:double-point-generic} Let
$F$ be a meromorphic function of
two variables, not singular at the origin, with
$F(z,w)=\numer(z,w)/\denom(z,w)=\sum_{r,s} a_{rs}z^rw^s$.

Suppose that $(z^*,w^*)$ is a strictly minimal, double point of $\sing$.  Let $\Hess(z^*,w^*)$ denote the
Hessian of $H$ at $(z^*,w^*)$ and suppose that $\det \Hess(z^*,w^*)\neq 0$.

Then 
\begin{enumerate} \romenumi
\item
for each \interior subset $K$ of $\cone (z^*,w^*)$, there is $c > 0$ such that

$$
a_{rs} = (z^*)^{-r}(w^*)^{-s} \left ( \frac{G(z^*,w^*)}{\sqrt{-(z^*)^2(w^*)^2 
   \det \Hess(z^*,w^*)}} + O(e^{-c |(r,s)|}) \right ) \qquad
\text{uniformly  for $(r,s)\in K$}.
$$

\item if $\delta = r/s$ lies on the boundary of $\overline{\cone(z^*,w^*)}$, then in direction $\delta$ there is a complete asymptotic expansion 
$$
a_{rs}  \sim (z^*)^{-r}(w^*)^{-s} \sum_{k\geq 0} b_k s^{-(k+1)/2}
$$
where $b_0 = \frac{G(z^*,w^*)}{2\sqrt{-(z^*)^2(w^*)^2 \det \Hess(z^*,w^*)}}.$

\end{enumerate} 
\romenumi

\end{theorem}

The geometric significance of the nonsingularity of $\Hess(z^*,w^*)$ is that
this is equivalent to the curves $\sing_0, \sing_1$ intersecting transversally.
If, on the other hand, $\det \Hess(z^*,w^*)=0$ (equivalently the curves are
tangent), then higher order information is required in order to compute
the relevant asymptotics. We treat the latter situation in Theorem~\ref{thm:degenerate}.

Note that in the situation of the above theorem, the asymptotic
exponential rate of $a_{rs}$ is constant on the interior of the cone
$\cone(z^*,w^*)$, that is, $a_{rs} \sim \exp (- (r,s) \cdot \vv)$ where $\vv
= (\log z^* , \log w^*)$ is constant on $\cone (z^*,w^*)$; this  differs
considerably from the asymptotics previously derived  for smooth
points~\cite{pemantle-wilson-multivariate1}.

\begin{example}[a lemniscate] \tlabel{eg:figure 8} Consider $F=1/\denom$, 
where $\denom(z,w)=19-20z-20w+5z^2+14zw+5w^2-2z^2w-2zw^2+z^2w^2$. The real points of $\sing$ are shown in Figure~\ref{fig:eight}.
\begin{figure}
\centerline{\psfig{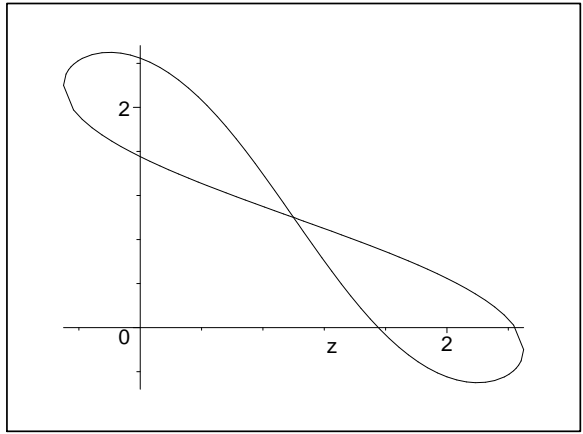}}
\vspace{10pt}
\caption{Figure-eight shape}
\label{fig:eight}
\end{figure}

The only common zero of $\denom$ and $\grad\denom$ is  
$(1,1)$, and so this is the only candidate for a double point. All other
 strictly
minimal elements of $\sing$  must be smooth. The second order part of the
 Taylor expansion of $H$ near $(1,1)$ is
$4(w-1)^2+10(z-1)(w-1)+4(z-1)^2=4[w-(3-2z)][w-(3/2-z/2)]$. The Hessian
determinant at $(1,1)$ is therefore $-36$.  The two tangent lines have slopes
$-1/2$ and $-2$, so that asymptotics in directions $r/s=\kappa\in [1/2,2]$
 cannot
be obtained by smooth points. By Theorem~\ref{thm:double-point-generic},
for slopes in the interior of this interval, the asymptotic $a_{rs}\sim 1/6$ 
holds, whereas on the boundary of the cone (slopes $1/2$ or $2$),  
$a_{rs}\sim 1/12$.
Directions corresponding to slopes outside the interval $[1/2, 2]$ may be
treated by the methods of \cite{pemantle-wilson-multivariate1}.  The 
exponential order is nonzero in these directions.

This example demonstrates that the local behavior of $H$ at 
the strictly minimal point $(1,1)$ determines asymptotics. 
In fact $H(1+z, 1+w)$ is homogeneous and irreducible in $\mathbb{C}[z,w]$ 
and hence does not factor globally (in the power series
ring $\mathbb{C}[[z,w]]$). 

\noproof
\end{example}

The previous theorem treated the simplest case, when $d=n=1$, and our 
subsequent results will be labelled in a similar way.  We begin with some
theorems involving nondegenerate points. The simplest case is when the
point is completely nondegenerate (a generalization of the hypothesis of 
Theorem~\ref{thm:double-point-generic}). 

\begin{theorem}[$d=n$; completely nondegenerate]\tlabel{thm:compnondeg}
Let $F$ be a function of $d+1$ variables  with $F(\b{z}) =
\numer(\b{z})/\denom(\b{z}) = \sum_{\b{r}}a_{\b{r}}\b{z}^\b{r}$.  
Suppose that $\b{z^*}$ is a strictly minimal, completely nondegenerate
multiple point of $\sing$ and that $F$ is meromorphic in a neighborhood
of $\torus (\b{z^*})$.  

Then there is $C$ such that for each \interior subset $K$ of
$\cone(\b{z^*})$, there is $c > 0$ satisfying 
$$
a_{\b{r}}\sim (\b{z^*})^{-\b{r}} \left ( C + O(e^{- c |\rr|}) \right )
\text{ uniformly for  $\rr \in K$}. 
$$
The constant $C$ is given by
$$
C = \frac{d! \, \phi(\b{z^*})}{|\det \cmat (\b{z^*})|}.
$$

\end{theorem}

\begin{unremark}
The picture to  keep in mind here is of $d+1$ complex hypersurfaces in
$\CC^{d+1}$ intersecting at a single point.  Notice that we again
obtain asymptotic constancy on the interior of the cone of allowable
directions, provided $\numer(\b{z^*})\neq 0$. If $\b{r}$ belongs to the boundary
 of $\cone$,  many different  asymptotics are possible.
\noproof 
\end{unremark}

In case $d<n$, a modification of the last result is required,
involving some linear algebra which we now introduce. 
Let $\simp=\simp_n$ denote the standard $n$-simplex in
$\mathbb{R}^{n+1}$,
$$
\simp := \left \{ \bb{\alpha} \in \mathbb{R}^{n+1} \left |
   \sum_{j=0}^{n}\alpha_j = 1 \right. \text{ and } \alpha_j \geq 0 
   \text{ for all } j \right \}.
$$
The subspace $\b{1}^\perp$ parallel to $\simp$ will 
be denoted by $\hyper$.
Assuming nondegeneracy of a minimal point $\b{z^*}$, the rows of $\cmat(\b{z^*})$
are $n+1$ vectors spanning $\R^{d+1}$, whence for any $\direc\in
\mathbb{R}^{d+1}$  the space of solutions  
$\mathcal{A}:= \mathcal{A}(\direc) := \{
\bb{\alpha} \mid \bb{\alpha} \cmat (\b{z^*})  = \direc \}$ will always be
$(n-d)$-dimensional.   The set $\mathcal{A} \cap \simp$ is just
the set of coefficients of convex combinations of rows of $\cmat$ that
yield the given vector $\direc$.  If $\direc \in \cone_0$ then
$\mathcal{A}$ is a subset of the translate of $\hyper$
containing $\simp$, thus if $\direc$ is in the interior of $\cone_0$ then
$\mathcal{A} \cap \simp$ is $(n-d)$-dimensional.  Let
$\mathcal{A}^\perp$  denote the orthogonal complement in $\hyper$ of
$\mathcal{A}_0$, the set $\mathcal{A}$ translated to the origin.  The
dimension (later denoted $\rho$) of
$\mathcal{A}^\perp$ is $d$.  The space  $\mathcal{A}^\perp$ does not
depend on $\direc$, since changing $\direc$ only translates
$\mathcal{A}$.  

The role that $\mathcal{A}$ will play is this: there will be an integral
over $\simp$ whose stationary points are essentially the points of
$\mathcal{A}$, each contributing the same amount; thus the total 
contribution will be the measure of $\mathcal{A}$ in the  
{\em complementary measure} which we  denote by $\sigma$.  The proper
normalization in the  definition of $\sigma$ is
$$\sigma (S) := \sigma (S \cap \mathcal{A} \cap \simp)
   := \frac{\mu_{n-d} (S \cap \mathcal{A} \cap \simp)} {\mu_n (\simp)}
$$
where $\mu_k$ denotes $k$-dimensional volume in $\mathbb{R}^{n+1}$.  One may
think of this as: $\sigma \times \mu_d =$ the normalized volume
measure, $\mu$, on $\simp$ ({\em Warning}: when $\mathcal{A}$ is a single point, it
is tempting to assume that $\sigma (\mathcal{A}) = 1$, but in
fact in this case $\sigma (\mathcal{A}) = \mu_n(\simp)^{-1}$.)   

\begin{proposition} \tlabel{prop:perp dim}
Assume that the columns of $\cmat$ are linearly independent.  Then
the projection of the column space of $\cmat$ onto $\mathcal{A}^\perp$
has dimension $d$.  
\end{proposition}

\begin{proof}
The matrix $\cmat$ is the sum of three columnwise projections, namely
one onto $\mathcal{A}_0$, one onto $\mathcal{A}^\perp$ and 
one onto the span of $\b{1}$.  The sum of the first two projections
annihilates the last column of $\cmat$, mapping the column space
to a space of dimension $d$.  By definition of $\mathcal{A}$, the space
$\mathcal{A} \cmat$ is a single point, hence $\mathcal{A}_0 \cmat = 0$, 
meaning that the first projection is null.  Therefore, the second 
projection has rank $d$.  
\end{proof}

\begin{definition} \tlabel{defn:cmat}
When the columns of $\cmat$ are linearly independent, let
$\overline{\cmat}$ denote the matrix representing the  linear
transformation $\vv \mapsto \vv \cmat$ on $\mathcal{A}^\perp$ with
respect to some basis of $\mathcal{A}^\perp$.  Note that
$\overline{\cmat}$ is independent of $\direc$ since $\mathcal{A}^\perp$
is, and that $\det \overline{\cmat}$ is independent of the choice of basis
of  $\mathcal{A}^\perp$.  For general $\cmat$, we extend this definition
so that $\overline{\cmat}$ is the projection of the column space of
$\cmat$ onto the space $\mathcal{A}^\perp$, then represented in a (fixed
but arbitrary) orthonormal basis of $\mathcal{A}^\perp$.  
\noproof
\end{definition}

In stating subsequent results, we shall rely increasingly on derived 
data such as $\overline{\cmat}$.  It is possible in principle  to 
give formulae for these asymptotics in terms of the original data 
$\numer$ and $\denom$, but such expressions rapidly become too 
cumbersome to be useful.  

\begin{theorem}[$\rho=d\leq n$; nondegenerate] \tlabel{thm:nondeg}
Let $F$ be a function of $d+1$ variables, with
$F(\b{z})=\numer(\b{z})/\denom(\b{z})=\sum_{\b{r}}a_{\b{r}}\b{z}^\b{r}$. 
Suppose that $\b{z^*}$ is a strictly minimal, nondegenerate multiple point
of degree $n+1$ of $\sing$ and that $F$ is meromorphic in a  neighborhood
of $\torus (\b{z^*})$.  

Then $\cone(\b{z^*})$ is a finite union of cones $\cone_j$ such that for each
\interior subset $K$ of each $\cone_j$, there are $c>0$ and a polynomial $P$
of degree at most $n-d$, such that
$$
a_{\b{r}} = (\b{z^*})^{-\b{r}} \left ( P(\rr) + O(e^{-c |\rr|}) \right )
\text{\qquad uniformly for $\rr \in K$.} 
$$
Here
$$
P(\rr) = \frac{\phi(\b{z^*})\sigma(\mathcal{A}(\direc) \cap \simp)}{|\det
\overline{\cmat(\b{z^*})}|} \, (r_{d+1})^{n-d}  + O((r_{d+1})^{n-d-1})\, , 
$$
where $\mathcal{A}(\direc)$ is the solution set to $\bb{\alpha} \cmat(\b{z^*})
= \direc$ with $\direc = \rr / r_{d+1} \in \cone_0(\b{z^*})$, and $\sigma$ 
is the complementary measure.

\end{theorem}

\begin{unremark} 

The approximation to $a_\rr$ is actually
piecewise polynomial and is asymptotically valid
throughout the interior of $\cone$.  The correction term may,
however, fail to be exponentially small on some lower-dimensional
surfaces in the interior of $\cone$ where the piecewise polynomial
is pieced together.  

We have broken the coordinate
symmetry in the formula for the leading term by parametrizing $\rr$
in terms of $r_{d+1}$ and $\direc$. See section~\ref{sec:future} for more comments.
\noproof
\end{unremark}

\begin{example}[3 curves in $2$-space]\tlabel{eg:distinct-tangents}
The simplest interesting  case to which Theorem~\ref{thm:nondeg}
applies is that where $\sing$ is the union of $3$ curves in $2$-space
which intersect at a strictly minimal point. 
Suppose that the strictly minimal singularity in question is at
$(1,1)$. In this case since $d=1$ and $n=2$, Theorem~\ref{thm:nondeg}
shows that 
$$
a_{rs}\sim  P(r,s)
$$ 
for some piecewise polynomial $P$ of degree at most $1$, whenever 
$\delta:= r/s$ lies in the
interior of the interval formed by the slopes $c_j=v_j'(1)$. We shall
obtain a more explicit expression for $P$.

Without loss of generality, we suppose that
$0 < c_0 \leq c_1 \leq c_2$, at least one of the two  inequalities being
strict.  Assume for now that  both inequalities are strict:
$c_0 < c_1 < c_2$.   Let
$\mathcal{A}_\delta=\{(\alpha_0,\alpha_1,\alpha_2)\in \simp \mid \alpha_0
c_0 + \alpha_1 c_1 + \alpha_2 c_2 = \delta\}$. When $\delta$ belongs to the
convex hull of $c_0, c_1$ and $c_2$, the affine set $\mathcal{A}_\delta$
is a line segment whose endpoints are on the boundary of the simplex
$\simp$.  Since $\mathcal{A}$ is orthogonal to $\b{c}$, the set
$\mathcal{A}^\perp$ is parallel to the projection $\bar{\b{c}}$ of
$\b{c}$ onto $\hyper$.   Letting $\mu$ denote the mean of the $c_j$,
$\Sigma/\sqrt{3}$  the standard deviation, we see that
the projection  is $\b{c} - \mu \b{1}$ and its euclidean length is $\Sigma$.

If $c_0<\delta<c_1$, then one
endpoint of the line segment is on the face $\alpha_2=0$, with
$\alpha_0=\frac{c_1-\delta}{c_1-c_0},\alpha_1=\frac{\delta-c_0}{c_1-c_0}$.
The other endpoint is when $\alpha_1=0$, with
$\alpha_0=\frac{c_2-\delta}{c_2-c_0},
\alpha_2=\frac{\delta-c_0}{c_2-c_0}$. The squared euclidean length of
this line segment simplifies to
$$3\frac{(\delta-c_0)^2 \Sigma^2}{(c_1-c_0)^2 (c_2-c_0)^2}.
$$
A similar argument gives the answer when $c_1<\delta <c_2$. The squared
euclidean length of the line segment is then
$$3\frac{(\delta-c_2)^2 \Sigma^2}{(c_2-c_1)^2 (c_2-c_0)^2} .  
$$
Both answers agree, and by continuity give the correct answer, when
$\delta=c_1$.  Since the area of $\simp_2$ is $\sqrt{3}/2$, 
the complementary measure $\sigma$ of $\mathcal{A}_\delta$ is 
$2/\sqrt{3}$ times its euclidean length. Thus we obtain
$P(r,s) = sf(r/s)$, where
$$
f(\delta) = 
\begin{cases}
2 \frac{\delta - c_0}{(c_1 - c_0) (c_2 - c_0)}, \qquad \text{if 
$c_0\leq \delta\leq c_1$};\\
2 \frac{c_2 - \delta}{(c_2 - c_1) (c_2 - c_0)}, \qquad \text{if 
$c_1\leq \delta\leq c_2$}.
\end{cases}
$$ 
Thus 
\begin{equation}
\label{eq:P}
P(r,s) = 
\begin{cases}
2 \frac{r - c_0 s}{(c_1 - c_0) (c_2 - c_0)}, \qquad 
\text{if $c_0\leq r/s \leq c_1$};\\
2 \frac{c_2 s - r}{(c_2 - c_1) (c_2 - c_0)}, \qquad 
\text{if $c_1\leq r/s \leq c_2$}.
\end{cases}
\end{equation}
Note that $f=0$ on the boundary of the cone, namely when $\delta=c_0$ or
$\delta=c_2$. The above formula extends to the case $c_0=c_1$ or
$c_1=c_2$ in the obvious way. This example illustrates the piecewise
polynomial nature of the asymptotics in $\cone$.  The approximation 
$a_{rs} \sim P(r,s)$ is exponentially close for $\delta$ in compact
subintervals of $(c_0 , c_1) \cup (c_1 , c_2)$, but requires polynomial
correction when $\delta \to c_1$ as well as at the endpoints $c_0$ and
$c_2$. 
\noproof  
\end{example}

More general formulae for $P$ in this case are described in the forthcoming work~\cite{baryshnikov-pemantle}. 

We move on to investigate situations which are degenerate according to our terminology.  We begin with the simplest case, namely when the minimal point is transverse.  Note
that in this case, the set $\mathcal{A}$ of linear combinations of the rows of $\cmat$ that yield
$\direc$ is always a single point since the rows of $\cmat$ are linearly independent.

\begin{definition} \tlabel{def:M}
If $\b{z^*}$ is a strictly minimal multiple point and $\bb{\alpha}$
is an element of $\simp$, we define the matrix $Q = Q(\b{z^*} , \bb{\alpha})$
to be the Hessian matrix of the function
$$
\bb{\w{\theta}} \mapsto - \log \bb{\alpha} \b{v}  (\b{\w{z^*}} e^{i \bb{\w{\theta}}} ) \, .
$$
at the point $\bb{\w{\theta}} = \b{0}$.  We define the matrix $M$ by
\begin{equation}\tlabel{eq:big hessian}
M := M(\b{z^*} , \bb{\alpha}) := \left(\begin{matrix}0&-i \overline{\cmat}(\b{z^*}) \\ 
   -i \overline{\cmat}^T(\b{z^*})&Q(\b{z^*}, \bb{\alpha}) \end{matrix}\right).
\end{equation}
\noproof
\end{definition}

\begin{theorem}[$\rho=n\leq d$; transverse] \tlabel{thm:transverse} 
Let $F$ be a function of $d+1$ variables, with 
$F(\b{z})=\numer(\b{z}) / \denom(\b{z})=\sum_{\b{r}}a_{\b{r}}\b{z}^\b{r}$. 
Suppose that $\b{z^*}$ is a strictly minimal, transverse multiple point
of degree $n+1$ of $\sing$ and that $F$ is meromorphic in a neighborhood
of $\torus(\b{z^*})$. 

Let $K$ be a  \interior subset of
$\cone_0 (\b{z^*})$ such that $\det M(\b{z^*} , \bb{\alpha (\direc)})\neq 0$ 
for $\direc \in K$, where $\bb{\alpha} (\direc)$ 
is the unique point in $\mathcal{A} (\direc)$. 
Then there is a complete asymptotic expansion 
$$
a_{\b{r}}\sim (\b{z^*})^{-\b{r}} \sum_{k\geq 0} b_k
(r_{d+1})^{\frac{n-d}{2}-k} \; \; \text{ uniformly for  $\direc \in K$}.
$$

If $\numer(\b{z^*})\neq 0$, the leading coefficient is given by
$$
b_0 = \frac{n!} {\sqrt{n+1}} \, \frac{(2\pi)^{\frac{n-d}{2}} \phi(\b{z^*})}
{\det M(\b{z^*}, \bb{\alpha (\direc)})^{1/2}}
$$
where the square root is the product of principal square roots of the eigenvalues.
\end{theorem}

\begin{example}[2 planes in 3-space]
Consider the trivariate sequence $(a_{r,s,t})$ whose terms are zero if any index 
is negative, and is otherwise given by the boundary condition 
$a_{0,0,0} = 1$ and the recurrence
\begin{eqnarray*}
16 a_{r,s,t} & = &  12 a_{r-1,s,t} + 12 a_{r,s-1,t} + 8 a_{r,s,t-1}
   -5 a_{r-1,s-1,t} - 3 a_{r-1,s,t-1} -3 a_{r,s-1,t-1} \\
&& -2 a_{r-2,s,t} - 2 a_{r,s-2,t} - a_{r,s,t-2}. 
\end{eqnarray*}

Let $F(x,y,z)=\sum_{r,s,t} a_{r,s,t} x^r y^s z^t$ be the associated 
generating function. Then we obtain directly the explicit form
$$
F(x,y,z)=\frac{16}{(4-2x-y-z)(4-x-2y-z)}.
$$
The point $(1,1,1)$ is a multiple point of $\sing$ which is readily seen
to be strictly minimal.  In fact this point is on the line segment
of multiple points $\{ (1,1,1) + \lambda (-1,-1,3) : -1/3 < 
\lambda < 1 \}$.  We focus here on the point $(1,1,1)$ in particular
because it is the only point giving asymptotics which do not decay 
exponentially.  Near $(1,1,1)$, the set $\sing$ is parametrized by
$z = 4 - 2 x - y$ and $z = 4 - x - 2 y$.  The cone $\cone$ is 
the positive hull of $(2,1,1)$ and $(1,2,1)$, with $\cone_0$ 
being the line segment between these.  Other multiple points
on the line segment govern other two dimensional cones.  For 
example, the point $(1/3, 1/3, 3)$ corresponding to $\lambda = 2/3$
governs the cone $\cone_{2/3}$ which is positive hull of 
$(2/3 , 1/3 , 3)$ and $(1/3 , 2/3 , 3)$, or in other words,
$\cone_{2/3}$ is the set of $(r,s,t)$ with $r+s = t/3$ and
$r , s \geq t/9$.  These other cones sweep out, as $\lambda$
varies, a polyhedral cone with non-empty interior.  There are
still directions outside the cone, namely directions $(r,s,t)$
with $\min \{ r , s \} \leq (r+s)/3$; asymptotics in these
directions are governed by smooth points on one of the 
two hyperplanes.  We return to the examination of the
non-exponentially decaying asymptotic directions.

The multiple point $(1,1,1)$ yields asymptotics for $a_{rst}$ with
$r+s = 3t$ and $r,s \geq t$.  Given $\direc =  (2 - \alpha , 1 + \alpha ,
1) \in \cone_0$, the set $\mathcal{A}$ is  the single point $(1 - \alpha
, \alpha)$.  The complement $\mathcal{A}^\perp$ is always all of
$\hyper$, which has an orthonormal basis $\{\xx\}$,  $\xx := (\sqrt{1/2}
, -\sqrt{1/2})$.

The matrix $\overline{\cmat}$ is equal to $((2,1) \cdot \xx , (1,2) \cdot
\xx) = \xx$.  This leads to $\det M = (A+B+C+D)/2$, where $A,B,C$ and $D$
are the entries of the restricted Hessian, $Q$.  A routine computation
shows that $A+B+C+D = 12$, leading to $\det M = 6$ (as is always the case
for transverse intersections of the pole set, the determinant of
$M$ does not vary with direction inside the cone of a fixed multiple
point).   Computing $\phi (\b{1})$ 
from~(\ref{eq:comp phi}) then gives 
$$
\phi (\b{1}) = \frac{2!} {(-1)^2} \, \frac{16} {2} = 16 \, .
$$
Plugging this into Theorem~\ref{thm:transverse} we find 
(to first order, as $t\to \infty$) that 
$$
a_{r,s,t} \sim \frac{16}{\sqrt{24\pi t}} \qquad \text{ if $r+s=3t$ and
$r/s\in (1/2, 2)$}.
$$ 
The above first order approximation differs from the true value of $a_{rst}$ by less than $0.3\%$ when $(r,s,t) = (90,90,60)$. 
\noproof
\end{example}

As mentioned in the introduction, instead of producing enough variants of
these theorems for a complete taxonomy, we will stop with just one  more
result.  The case we describe is the most degenerate, namely  where all
the sheets $\sing_j$ are tangent at $\b{z^*}$, so that $\cone(\b{z^*})$ is a
single ray.  In this case $\det M$ varies with $\bb{\alpha}$ and the resulting formula is an integral over $\bb{\alpha}$.  The methods used here may be adapted to prove results when the degeneracy is somewhere between this extreme and the nondegenerate cases.  

\begin{theorem}[$\rho=0$] \tlabel{thm:degenerate}
Let $F$  be a function of $d+1$ variables, not singular at the origin, with
$F(\b{z})=\numer(\b{z})/\denom(\b{z})=\sum_\b{r} a_\b{r}\b{z}^{\b{r}}$. 
Suppose that $\b{z^*}$ is a strictly minimal, multiple point of 
degree $n+1$ of $\sing$ and that $F$ is meromorphic in a neighborhood
of $\torus (\b{z^*})$. Further suppose that $\cone(\b{z^*})$ 
is a single ray.  

If $\det Q(\b{z^*}, \bb{\alpha})\neq 0$ on $\simp$, then for $\b{r}\in \cone(\b{z^*})$
we have 

$$
a_\b{r} = \b{z^*}^{-\b{r}} [b_0 (r_{d+1})^{n-d/2}  + 
O(r_{d+1})^{n - d/2 - 1}] \, .
$$

The value of $b_0$ is given by 
$$
b_0 = \frac{\phi(\b{z^*})}{(2\pi)^{d/2}}\int_\simp \det
Q(\b{z^*}, \bb{\alpha})^{-1/2}\,d\meas(\bb{\alpha}).
$$
\end{theorem}

\begin{unremark}
Theorem~\ref{thm:degenerate} covers the situation where all the sheets coincide, as opposed to merely being tangent. That situation can also be analyzed  in a different way by a slight modification of the proof of  the smooth case, as discussed in~\cite{pemantle-wilson-multivariate1}. The methods of~\cite{pemantle-wilson-multivariate1} yield a complete asymptotic expansion. 
\noproof
\end{unremark}

\begin{example}
\tlabel{eg:tangent2D}

The simplest case illustrating Theorem~\ref{thm:degenerate} is that of
two curves $w  = u_j(z)$ in $\mathbb{C}^2$, intersecting tangentially at the
strictly minimal point $\b{z}=(1,1)$. Suppose that $\re  \log
v_j(e^{i\theta})  =  -d_j \theta^2 + \cdots$ with $d_j >0$. The simplex in
question is one-dimensional and we identify it with the interval $0\leq
t\leq 1$. Then for $0\leq t\leq 1$ the matrix $M(\b{z}, t)$ is just the
$1\times 1$ matrix with entry $d_t=(1-t)d_0+td_1$. Hence when $d_0\neq d_1$ we
obtain asymptotics in the unique direction $\direc$ of the theorem: 

\begin{align*}
a_{rs}&\sim \frac{\phi(1,1)\sqrt{s}}{\sqrt{2\pi}} \int_0^1 (d_t)^{-1/2}\, dt\\
&= \frac{\phi(1,1)\sqrt{s}}{\sqrt{2\pi}} \frac{1}{d_1- d_0}\int_{d_0}^{d_1}
y^{-1/2} \,dy\\
&= \frac{2\phi(1,1)\sqrt{s}}{\sqrt{2\pi}} \frac{\sqrt{d_1}-\sqrt{d_0}}{d_1-
d_0}\\
&= \frac{\phi(1,1)\sqrt{s}}{\sqrt{2\pi}} \frac{2}{\sqrt{d_0}+\sqrt{d_1}}.
\end{align*}

This final formula is easily seen to hold also in the case $d_0=d_1$,
in which case the formula $\phi (1,1)\sqrt{s} / \sqrt{2 \pi d_0}$ agrees
with the formula in~\cite{pemantle-wilson-multivariate1}.  
\noproof
\end{example}

\begin{unremark}
Since analysis in $1$ variable is considerably easier than in general,
it is possible to derive a result for plane curves even when $d_j$
 vanishes.  We omit the details here, but 
note that such a derivation is possible because (by minimality)
 we must have $\log wv_j(
z e^{i\theta})=ic_j\theta -d \theta^m + \cdots$, where $c_j\geq 0$, $m$ is
even and $\re  d  > 0$.

If all tangents are equal at $(z^*,w^*)$, then $\cone_0 (z^*,w^*)$ is a singleton.
If all $\re  \log w v_j (z e^{i\theta}) $ vanish to the
same exact order $m$ and  $(r,s) \in \cone(z^*,w^*)$, then
$$
a_{rs} = (z^*)^{-r}(w^*)^{-s} [b_0 s^{n-1/m} + O(s^{n - 2/m})]
$$
with $b_0$ given by 
$$
b_0=\frac{\phi(z^*,w^*) \Gamma(1/m)}{2\pi(1-1/m)}
\frac{d_1^{1-1/m}-d_0^{1-1/m}}{d_1-d_0}.
$$
The above formula for $b_0$ also holds in the limit when $d_0 = d_1$, when
it becomes $\displaystyle{b_0=\frac{\phi(z,w)\Gamma(1/m)}{2\pi(d_0)^{1/m}}}$.
\noproof
\end{unremark}

\section{Proofs} \tlabel{sec:proofs}
Throughout, we assume that $\b{z^*}$ is a fixed strictly minimal
multiple point of $\sing$, of multiplicity $n+1$, and $F(\b{z}) =
\phi(\b{z}) / \prod_{j=0}^{n} (1 - z v_j(\w{\b{z}}))$ near $\b{z^*}$ as
in Definition~\ref{defn:preparation}. Recall our notational conventions
$s := r_{d+1}, z := z_{d+1}$. To ease notation we shall
later assume that $\b{z^*} = \b{1}$.

The proofs follow a $4$-step process. First we use the residue theorem
in one variable to reduce the dimension by $1$ and to restrict attention
to a neighbourhood of $\b{\w{z^*}}$. Second, the residue sum resulting
is rewritten in a form more amenable to analysis. Third, the Cauchy-type
integral is recast in the Fourier-Laplace framework by the usual
substitution $\b{\w{z}} = \b{\w{z^*}} e^{i\bb{\theta}}$. Finally, detailed
analysis of these Fourier-Laplace integrals (including some
generalizations of known results we have included separately in 
\cite{pemantle-wilson-analysis} to save space here)  yields our desired
results.

\subsection*{Restriction to a neighborhood of $\b{z^*}$}

We first apply the residue theorem in one variable: the proof of the following
proposition is entirely analogous to that of
\cite[Lemma~4.1]{pemantle-wilson-multivariate1}. 

\begin{proposition}[Local residue formula]
\tlabel{prop:localize} 
Let $R(s,\b{\w{z}},\varepsilon)$ denote the (finite) sum of the residues of 
$g: z\mapsto z^{-s-1} F(\b{\w{z}}, z)$ at $z = u_j(\b{\w{z}})$ inside the
ball $|z - z^*| < \varepsilon$. 

Let $\nbd'$ be a neighborhood of $\w{\b{z^*}}$ in $\torus(\w{\b{z^*}})$.
Then for sufficiently small $\varepsilon > 0$,
\begin{equation} \tlabel{eq:tores}
|\b{z}^{\rr}| \left|a_\rr - (2\pi i)^{-d} \int_{\nbd'} - 
\b{ \w{z} }^{-\b{ \w{r} } - \b{1} }  R(s,\b{\w{z}},\varepsilon) 
 \, d\b{ \w{z} } \right |   
\text{ is exponentially decreasing in $s$ as $s \to \infty$.}
\end{equation}
This is uniform in $\b{r}$ as $\b{r}/s$ varies over some neighborhood of 
$\cone_0(\b{z^*})$.
\noproof
\end{proposition}

\subsection*{Rewriting the residue sum}

After using the local residue formula \eqref{eq:tores} we must  perform a
$d$-dimensional integration of the residue sum $R$. The following lemma, whose
proof may  be found in~\cite[p. 121, Eqs 7.7 \&
7.12]{devore-lorentz-constructive}, will yield a more tractable form for $R$. 
Recall from Section~\ref{sec:results} the relevant facts about
simplices. In the following, the notation $\bb{\alpha} \b{v}$ denotes
the scalar product $\sum_{j=0}^n \alpha_j v_j$, and $h^{(n)}$ the $n$th
derivative of $h$.

\begin{lemma} \tlabel{lem:simplex} Let $h$ be a function of one complex
variable, analytic at $0$, and let $\meas$ be the normalized volume
measure on $\simp_n$. Then 
$$
\sum_{j=0}^{n} \frac{h(v_j)}{\prod_{r \neq j} (v_j - v_r)} =
\int_{\simp_n} h^{(n)} (\bb{\alpha}\b{v}) \, d\meas (\bb{\alpha}) 
$$
both as formal power series in $n+1$ variables $v_0, \dots ,v_{n}$ and in a 
neighborhood of the origin in $\mathbb{C}^{n+1}$. 
\noproof
\end{lemma} 

\begin{corollary} [Residue sum formula]
\tlabel{cor:residue sum}
Let $R(s,\w{\b{z}},\varepsilon)$ be
the sum of the residues of the function $g: z\mapsto z^{-s-1}F(\b{z})$
inside the ball $|z - z^*| < \varepsilon$. Define $h_{s,\w{\b{z}}} (y) =
y^{s+n} \phi(\w{\b{z}}; 1/y)$.

Then for sufficiently small $\varepsilon$, there is $\delta > 0$ such
that for $|\w{\b{z}} - \w{\b{z^*}}| < \delta$, 
\begin{equation}
\tlabel{eq:residue} R(s,\w{\b{z}},\varepsilon) = \int_{\simp}
h_{s,\w{\b{z}}}^{(n)} (\bb{\alpha}\b{v})\,d\meas(\bb{\alpha}).
\end{equation}
\end{corollary}

\begin{proof} First suppose that the
functions $v_0, \dots, v_{n}$ are distinct.  Choose $\varepsilon$
sufficiently small and $\delta > 0$ such that $g$ has exactly $n+1$
simple poles in $ | z - z^* | < \varepsilon$ whenever 
$|\w{\b{z}} - \w{\b{z^*}}| < \delta$.  Then there are $n+1$ residues, the $j$th one being 
$$
\frac{v_j (\w{\b{z}})^{s+n} \phi (\w{\b{z}} , 1 / v_j
(\w{\b{z}}))}{\prod_{r \neq j}(v_r (\w{\b{z}}) - v_j (\w{\b{z}}))} = 
\frac{h_{s,\w{\b{z}}} (v_j(\w{\b{z}}))}{\prod_{r \neq j}  (v_r
(\w{\b{z}}) - v_j (\w{\b{z}}))}.
$$  
The result in this case now follows
by summing over $j$ and applying Lemma~\ref{lem:simplex}.    

In the case when $v_0, \dots , v_{n}$ are not distinct, let $v_j^t$ be
functions approaching $v_j$ as $t \rightarrow 0$, such that $v_j^t$ are
distinct  for $t$ in a punctured neighborhood of $0$, and let $g^t(z) =
z^{-s-1}\phi (\b{z}) / \prod_{j=0}^{n} (z - 1 / v_j^t (\w{\b{z}}))$, so
that $g^t \rightarrow g$ as well. The sum of the residues of $g^t$ may be
computed by  integrating $g^t$ around the circle $|z-z^*|=\varepsilon$, and
since $g^t \rightarrow g$, this sum approaches the  sum of the residues
of $g$.  Since the expression in~\eqref{eq:residue}  is continuous in the variables $v_j$, this proves the general case. 
\end{proof}

The residue sum formula and Proposition~\ref{prop:localize} combine to show that
for some neighborhood $\nbd'$ of $\b{\w{z^*}}$ in $\torus(\b{z^*})$, 
\begin{equation} \tlabel{eq:iterint}   
|\b{(z^*)^r}| \left| a_{\b{r}} - (2\pi i)^{-d} \int_{\nbd'} 
- \w{\b{z}}^{-\w{\b{r}} - \b{1}} \int_{\simp}h_{s,\w{\b{z}}}^{(n)}(\bb{\alpha}\b{v}) 
\,d\meas(\bb{\alpha})\,d\w{\b{z}}\right| \quad \text{ is exponentially 
decreasing in $s$},
\end{equation} 
where, as usual, we have let $s$ denote $r_{d+1}$.

\subsection*{Recasting the problem in the Fourier-Laplace framework}

In this subsection we assume throughout (purely in order to ease notational complexity)
that $\b{z^*} = \b{1}$.

\begin{lemma} [Reduction to Fourier-Laplace integral]
\tlabel{lem:tooscint}
For $0 \leq k \leq n$, $e^{i\w{\bb{\theta}}} \in \nbd'$ and 
$\bb{\alpha} \in \mathbb{R}^{n+1}$, define
\begin{eqnarray*}
p_k (s) & := & \frac{n! (s+n)!} {k! (n-k)! (s+k)!} \\[1ex]
f (\w{\bb{\theta}} , \bb{\alpha}) & := & \frac{i \w{\b{r}}} {s} 
   - \log (\bb{\alpha} \vv (e^{i \w{\bb{\theta}}}) ) \\[1ex]
\mult_k (\w{\bb{\theta}} , \bb{\alpha}) & := & \left. \left ( \frac{d} {dy} \right )^k 
 \phi (e^{i\w{\bb{\theta}}} , 1/y) \right|_{y = \vv (e^{i \bb{\w{\theta}}})} \\[1ex]
I(s; f, \mult_k) & := & \int_\nice e^{-s f}\mult_k \, d\mu_{n+d}.
\end{eqnarray*}
Then there is a neighborhood $\nbd$ of $\b{0}$ in $\mathbb{R}^d$, a product of 
compact intervals, such that
\begin{equation}
\tlabel{eq:multint}
\left| a_{\b{r}} - (2\pi)^{-d} \frac{n!}{\sqrt{n+1}} 
\sum_{k=0}^{n}p_k(s) I(s; f, \mult_k) \right| 
\quad\text{is exponentially decreasing in $s$}\, ,
\end{equation}
where $\nice := \nbd \times \simp_n$. 
\end{lemma}

\begin{unremark}
Note for later that each $\psi_k$ is independent of $\bb{\alpha}$, as is
 $f(\b{0}, \bb{\alpha})$.  
\end{unremark}

\begin{proof}
Recall that $\sigma (\mathcal{A}) = n ! / \sqrt{ n + 1 } = V_n(\simp_n)^{-1}$.
Applying Leibniz' rule for differentiating $h$ yields

\begin{align*}
h_{s, \w{\b{z}}}^{(n)}(y)&=\sum_{k=0}^n \binom{n}{k} \left(
   \frac{d} {dy} \right )^{n-k} y^{s+n}
   \left( \frac{d} {dy} \right )^{k} \phi(\w{\b{z}},1/y)\\
&=\sum_{k=0}^n \frac{n!(s+n)!}{k!(n-k)!(s+k)!} y^{s+k} 
   \left ( \frac{d} {dy} \right )^{k} \phi(\w{\b{z}},1/y)\\
&=y^s\sum_{k=0}^n p_k(s) y^k \left ( \frac{d} {dy} \right )^k
   \phi (\b{\w{z}} , 1/y) \, .
\end{align*}

Substituting this into~(\ref{eq:iterint}) with $\w{\b{z}}
= e^{i\w{\bb{\theta}}}$, shrinking $\nbd'$ if necessary and letting 
$\nbd$ be the inverse image under $\w{\bb{\theta}} \mapsto 
e^{i\w{\bb{\theta}}}$ completes the proof and the third step.  
\end{proof}

\subsection*{Analysis of the Fourier-Laplace integral} 

The proofs of all the theorems in Section~\ref{sec:results} proceed in parallel while we compute the asymptotics for~(\ref{eq:multint}) with the aid of the complex Fourier-Laplace expansion. To that end, observe that $\nice$ is a compact product of simplices and intervals of dimension $n+d$ inside $\R^{n+d+1}$, $f$ and $\psi_k$ are complex analytic functions on a neighbourhood of $\nice$ in $\hyper \times \R^d$ (this requires a sufficiently small choice of $\varepsilon$ in Proposition~\ref{prop:localize}), and $\re f\geq 0$ on $\nice$. 

To apply the standard theory, we must locate the critical points of $f$ on $\nice$. Recall from Definition~\ref{defn:cone} the matrix $\cmat = \cmat(\b{z^*})$ whose rows are extreme rays in the cone of outward normals to support hyperplanes to the logarithmic domain of convergence of $F$.

\begin{proposition}[Critical points of $f$]
\tlabel{prop:statset} For $\direc := \b{\w{r}} / s \in \cone_0$, 
let $\mathcal{A}(\direc)$ denote the solution set  $\{ \bb{\alpha}\in
\mathbb{R}^{n+1} \mid \bb{\alpha} \cmat = \direc \}$.  Let $\statset$
denote the set of critical points of $f$ on $\nice$.  Then $\statset$
consists precisely of those points  $(\b{0}, \bb{\alpha})$ with
$\bb{\alpha}\in \mathcal{A}(\direc) \cap \simp$. 
\end{proposition}

\begin{proof}
By strict minimality of $\b{z^*}$, for each $j$ the modulus of  $v_j
(\b{\w{z}})$ achieves its maximum only when $\b{\w{z}} = \b{\w{z^*}}$. 
Thus any convex combination of $v_j (\b{\w{z}})$ with $\b{\w{z}} \neq
\b{\w{z^*}}$ has modulus less than $|v_j (\b{\w{z^*}})|$.  In other words,
when $\bb{\w{\theta}} \neq \b{0}$, the real part of $f$ is strictly positive,
from which it follows that all critical points  of $f$ are of the form
$(\b{0} , \bb{\alpha})$ for $\bb{\alpha} \in \simp$.  

In fact $f$ is somewhat degenerate: $f (\b{0} , \bb{\alpha}) = 0$ for all
$\bb{\alpha} \in \simp$, so not only does the real part of $f$ vanish
when $\bb{\w{\theta}} = \b{0}$, but also the  $\bb{\alpha}$-gradient of
$f$ vanishes there.  We compute the  $\bb{\w{\theta}}$-derivatives at
$\bb{\w{\theta}} = \b{0}$ as follows.   For $1 \leq j \leq d$,
$$
\frac{\partial f}{\partial \theta_j} = i 
\left.
\left ( \frac{r_j} {s} - \frac{z_j}{z_{d+1}} \bb{\alpha}
\left ( \frac{\partial}{\partial z_j} \right )
\b{v}(\w{\bb{z}})\right )
\right |_{\b{\w{z}} = \b{\w{z^*}}} \, .
$$
Recalling the definition of $\cmat$, we see that these vanish 
simultaneously if and only if $\direc = \bb{\alpha} \cmat$, which
finishes the proof.
\end{proof}

The set of critical points of $f$ on $\nice$ is thus the intersection of $\nice$ with an affine subspace of 
$\mathbb{R}^{n+d+1}$. We now wish to use known results on asymptotics of 
Fourier-Laplace integrals. Note that in the situation of
Theorems~\ref{thm:double-point-generic},~\ref{thm:compnondeg}
and~\ref{thm:transverse}, $f$ has a unique critical point on $\nice$, whereas in Theorem~\ref{thm:nondeg}, $\statset$ is essentially an affine subspace of $\simp$, and in Theorem~\ref{thm:degenerate} $\statset$ is essentially just $\simp$. Though the existing literature  
on asymptotics of integrals is extensive, we have been unable to find in it results applicable to all cases of interest of us here. We have therefore derived the extra asymptotics  we need in \cite{pemantle-wilson-analysis} and cite the relevant ones below in Lemmas~\ref{lem:smallcrit} and \ref{lem:bigcrit}.

Recall that $\mathcal{A}^\perp$ is the orthogonal complement to (any and every) $\mathcal{A}(\direc)$ in $\hyper$; for each point $y=(\b{0}, \bb{\alpha})$ of $\statset$, let $B_y$ be the product of $\nbd$ with the affine subspace of $\simp$ containing $\bb{\alpha}$ and parallel to $\mathcal{A}^\perp$.

Then the integrals in Lemma~\ref{lem:tooscint} may be decomposed as 
\begin{equation} \tlabel{eq:decompose}
\int_\nice e^{-s f} \psi_k \, d\mu_{n+d} 
= \int_\statset \left(\int_{B_y} e^{-s f}  \psi_k \, d\mu_{\rho +d} \right) 
\, d\sigma(y) \, ,
\end{equation}
where $\rho := \dim \mathcal{A}^\perp$.

We require a more explicit description of the Hessian matrix of $f$.

\begin{proposition}[Hessian of $f$]
\tlabel{prop:hessian} 
Let $y= (\b{0}, \bb{\alpha}) \in \statset$.
At the point $y$, the Hessian of the restriction of $f$ to $B_y$ has the block form
\begin{equation}\tlabel{eq:hessian}
M(\b{z^*}, \bb{\alpha})=\left(\begin{matrix}0&-i \overline{\cmat(\b{z^*})}\\-i
\overline{\cmat(\b{z^*})}^T&Q(\b{z^*}, \bb{\alpha}) \end{matrix}\right).
\end{equation}

In this decomposition:
\begin{itemize}
\item the zero block  has dimensions $\rho\times \rho$, where
$\rho=\dim \mathcal{A}^\perp=\rk \overline{\cmat}(\b{z^*})$ 
\item the $j$th column of $\overline{\cmat}(\b{z^*})$ is the projection of the
$j$th column 
of $\cmat(\b{z^*})$ onto $\mathcal{A}^\perp$ (expressed in some orthonormal basis) 
\item $Q(\b{z^*}, \bb{\alpha})$ is the Hessian at $y$ of the restriction of $f$ to the
$\w{\bb{\theta}}$-directions, as defined in Definition~\ref{def:M}.
\end{itemize}
\end{proposition}

\begin{proof} 
Constancy of $f$ in the $\simp$ directions at $\bb{\w{\theta}} = \b{0}$
shows that the second partials in those directions vanish, giving 
the upper left block of zeros.  Computing $(\partial / \partial \theta_j)
f$ up to a constant gives 
$$
- \frac{i} {\bb{\alpha} \vv (\b{\w{z^*}})} \bb{\alpha} z_j \frac{\partial} 
   {\partial z_j} \vv (\b{\w{z^*}}) \, 
$$
and since $\bb{\alpha} \vv (\b{\w{z^*}})$ is constant when $\bb{\w{\theta}} 
= \b{0}$, differentiating in the $\bb{\alpha}$ directions recovers 
the blocks $-i \overline{\cmat}$ and $-i \overline{\cmat}^T$.   
The second partials in the $\bb{\w{\theta}}$ directions are of course
unchanged.  To see that the dimension of $\mathcal{A}^\perp$ is the
rank of $\overline{\cmat}$, observe that no nonzero element of 
$\mathcal{A}^\perp$ can be orthogonal to column space of 
$\cmat$, since then it would have been in $\mathcal{A}$; thus
the projection of the columns of $\cmat$ onto $\mathcal{A}^\perp$
spans $\mathcal{A}^\perp$.  
\end{proof}

We may now prove Theorems~\ref{thm:double-point-generic},~\ref{thm:compnondeg}
and~\ref{thm:transverse}, in reverse order from most general to most special. In
these cases the function $f$ in the representation~(\ref{eq:multint}) has only
one critical point, so the decomposition~(\ref{eq:decompose}) is not really
necessary. 

As mentioned above,  we have had to derive asymptotic expansions for some 
cases. We start with one of these, applicable to the situations of 
Theorems~\ref{thm:double-point-generic},~\ref{thm:compnondeg}, 
and~\ref{thm:transverse}. 

\begin{lemma}[Fourier-Laplace expansion for isolated critical point] 
\tlabel{lem:smallcrit}
Suppose that $f$ has a unique critical point $\b{0} \in \nice$, at which the
 Hessian $\Hess$ is nonsingular.

\begin{enumerate} \romenumi
\item
If $\b{0}$ has a neighborhood in $\nice$ diffeomorphic to $\R^d$  then  
as $s \to \infty$ there is an explicitly computable asymptotic expansion 
$$
I(s; f , \mult_k) \sim  \sum_{j=0}^\infty a_j s^{-d/2-j}. 
$$
In particular
$$
a_0 = \mult_k (\b{0}) \frac{(2 \pi)^{d/2}} {\sqrt{\det(\Hess)}}
$$
where the square root is defined to be the product of the 
principal square roots of the eigenvalues of $\Hess$.  
\item
If  $\b{0}$ has a neighborhood in $\nice$ diffeomorphic to a 
$d$-dimensional half-space  then there is an explicitly computable asymptotic
 expansion
$$
I(s ; f , \mult_k) \sim \sum_{j = 0}^\infty a_j s^{- d/2 - j/2}. 
$$
Here the value of $a_0$ is half that of the previous case.
\end{enumerate}

\end{lemma}

\begin{proof}
This is a straightforward corollary of \cite[Theorems 7.7.1, 7.7.5, 7.7.17 
(iii)]{hormander-partial-differential1}. Details are given in 
\cite{pemantle-wilson-analysis}.

\end{proof}

\subsection*{Proof of Theorem~\ref{thm:transverse}}
We prove the result in the case $\b{z^*} = \b{1}$. The result in the general case follows directly via the change of variable  $\b{z} \to \b{z}/ \b{z^*}$. 

The rank of $\cmat$ is $n+1$, so $\mathcal{A}(\direc)$ is a single point and
$f$ has a unique critical point, $(\b{0},\bb{\alpha})$. Hence $\mathcal{A}^\perp
= \hyper$ and has dimension $n$.  The critical point  is in the interior of
 $\nice$ as long as $\direc$ is in the interior of $\cone_0$, which occurs
 if and only if $\b{r}$ is in the interior of $\cone$.

Since $\mathcal{A}$ is a singleton, we may use the formula~(\ref{eq:multint})
 and Lemma~\ref{lem:smallcrit}. The Fourier-Laplace expansion of 
$\int_\nice e^{-s f} \psi_k$ will have terms
$s^{-(n + d)/2 - j}$ for all $j$, and when multiplied by the degree $n - k$
polynomial $p_k$, will contribute to the $s^{(n - d)/2 - k - j}$ terms for all
$j$.  Each of these terms is explicitly computable from the partial derivatives
of $\numer$ and $\denom$ via $\psi_k$, but these rapidly become messy.  The
leading term of $s^{(n - d)/2}$ comes only from $k = 0$: 
$$ a_\rr \sim (2 \pi)^{-d} \frac{n!} {\sqrt{n + 1}}
s^n \frac{(2 \pi)^{(n + d)/2} s^{(n - d)/2} \phi
(\b{z})} {\sqrt{\det M(\b{z^*}; \bb{\alpha})}} 
$$
which proves the theorem in the case $\b{z^*} = \b{1}$. 
\noproof 

\subsection*{Proof of Theorem~\ref{thm:compnondeg}}

Here $\rank (\cmat) = d+1 = n+1$ and so the proof of 
Theorem~\ref{thm:transverse}
applies. We can obtain a more explicit result in this case. The block
decomposition of $M(y)$ in~(\ref{eq:hessian}) is into four blocks of dimensions
$d \times d$, implying that $M(y)$ is nonsingular with determinant $(\det
\overline{\cmat})^2$ (in performing column operations to bring $M(y)$ to block
diagonal form, we pick up a factor of $(-1)^{d^2}$ which cancels the factor
$i^{2d}$).  The last column of $\cmat$ is $\b{1}$, whose projection onto
$\mathcal{A}^\perp$ is null; its projection onto the space spanned by $\b{1}$
 has length $\sqrt{d+1}$, whence $|\det{\cmat}| = \sqrt{d+1} \,|\det
 \overline{\cmat}|$.

Thus we obtain
$$
a_\rr = C  + O(s^{-1})
$$

where
$$
C = \frac{d! \, \phi (\b{z^*}) } {|\det \cmat(\b{z^*})|} \, .
$$
To finish the proof, we quote~\cite[Corollary~3.2]{pemantle-finite-dimension} to see
that in fact $a_\rr - C$ is exponentially small on \interior subcones of the interior of 
$\cone$.

\begin{unremark}
If $\direc$ is on the boundary of $\cone$, it may happen that
the point $(\b{0} , \bb{\alpha})$ has $\bb{\alpha}$ lying on
a face of $\simp$ of dimension $d - 1$.  In that case, using the
second formula from Lemma~\ref{lem:smallcrit} instead of the first  gives a leading term with exactly half the previous magnitude. Here we do not claim exponential decay of $a_\rr - C$. In Example~\ref{eg:figure 8}, this occurs when $r/s \in \{ 1/2 , 2 \}$.  
\end{unremark}

\subsection*{Proof of Theorem~\ref{thm:double-point-generic}}
The second case is covered by the previous remark.

Details for the first case follow. The analysis of Theorem~\ref{thm:compnondeg} applies. We  derive a more explicit result, first supposing that $(z^*, w^*) = (1,1)$.

For $j\in \{0,1\}$, write $v_j(e^{i\theta}) = i c_j\theta 
+ O(\theta^2)$,  where $c_j\geq 0$. The matrix $\cmat$ in this case is 
simply $(\begin{smallmatrix} c_0 & 1 \\ c_1 & 1 \end{smallmatrix})$.
It will be seen below that $\det \cmat = 0$ if and only if $\det \b{H} = 0$, so for the purposes of this theorem we may and shall  suppose that $c_0 \neq c_1$. 

We have 
$$
a_{rs} \sim \frac{\phi(1,1)} {|c_0 - c_1|} \qquad 
\text{uniformly for $r/s$ \interior to $[c_0, c_1]$}.
$$

This expression for the constant is quite handy, but to recover the
statement of the theorem, we also express it in  terms of the partial
derivatives of $\denom$.  Recall that  $\denom = \chi Q$ where $Q(z,w) =
(1 - w v_0(z))(1 - w v_1(z))$ and $\chi \neq 0$ near $(1,1)$.  Then at
$(1,1)$, $Q = Q_z = Q_w = 0$ and so $\denom_{zz} = \chi Q_{zz} = 2 \chi c_0
c_1$. Similarly, $\denom_{wz} = \chi(c_0 + c_1)$ and $\denom_{ww} =
2\chi$. Thus $c_0, c_1$ are the roots of the quadratic $(\denom_{ww}) x^2
- (2 \denom_{wz}) x + \denom_{zz} = 0$.  Solving via the quadratic
formula gives $|c_0 - c_1| = 2 \sqrt{- \det \Hess} / \denom_{ww}$.  This
shows that $\det \b{H} = 0$ if and only if $c_0 = c_1$, as asserted
above.

Formula~(\ref{eq:comp phi}) then yields
$$
a_{rs} \sim \frac{\phi (1,1)} {|c_0 - c_1|} = 
\frac{2 \numer (1,1)} {\denom_{ww} (1,1) |c_0 - c_1|} = 
\frac{\numer (1,1)} {\sqrt{ - \det \Hess (1,1)}} \, .
$$

The general result follows by a change of scale mapping the  multiple
point $(z^* , w^*)$ to $(1,1)$.  Details of the computation  are as
follows.  Make the change of variables $\widetilde{z} = z/z^*,
\widetilde{w} = w/w^*$, and write $\widetilde{G}$ for $G$ considered as a
function of $\widetilde{z}, \widetilde{w}$. In an obvious notation,
$\widetilde{a}_{rs} = (z^*)^r (w^*)^s a_{rs}$.  Then, evaluating at
$(\widetilde{z},\widetilde{w}) = (1,1)$, we obtain $H_{\widetilde{z}} =
z^*H_z$, $H_{\widetilde{w}} = w^*H_w$. Repeating this yields
$H_{\widetilde{z} \widetilde{z}} = (z^*)^2 H_{zz}$, 
$H_{\widetilde{w}{\widetilde{z}}} = z^*w^*H_{wz}$ and
$H_{\widetilde{w}{\widetilde{w}}} = (w^*)^2H_{ww}$.

The special case above now yields
$$
\widetilde{a}_{rs} \sim 
\frac{\widetilde{G}(1,1)} {\sqrt{ - \widetilde{D}(1,1)}} = 
\frac{G(z^*,w^*)} {\sqrt{ - (z^*)^2 (w^*)^2 \det \Hess(z^*,w^*)}}
$$
from which the result follows immediately.
\noproof

To prove Theorems~\ref{thm:nondeg} and~\ref{thm:degenerate} we require a result on asymptotics of integrals that allows for non-isolated critical points. The lemma below says that, as seems intuitively plausible, we may compute  asymptotic expansions in directions transverse to $\statset$ and then integrate their leading terms along $\statset$ to get the leading term of the original expansion.

\begin{lemma}[Fourier-Laplace expansion for continuum of critical points]
\tlabel{lem:bigcrit}
Let $y\in \statset$ and let $M(y)$ be the Hessian at $y$ of the restriction of $f$ to $B_y$.
If $\mult_k(y) \neq 0$ and $\det M(y)$ is nonzero on $\statset$ then
$$ 
\int_\nice e^{-s f} \mult_k  \, d\bb{\w{\theta}} \times d\mu
\sim \left(\frac{2\pi}{s}\right)^{(\rho +d)/2}  \int_\statset 
\frac{\mult_k (y)}{\sqrt{\det M(y)}} \,d\sigma(y). 
$$
\end{lemma}

\begin{proof}
This is a specialization of \cite[Thm 2]{pemantle-wilson-analysis} 
(the proof of that result is considerably more involved than in 
the isolated critical point case).  
\end{proof}

Lemma \ref{lem:bigcrit} yields the following expansion:
\begin{equation} \tlabel{eq:integrated}
a_\rr \sim (2 \pi)^{-d} p_0 (s) \left(\frac{2\pi}{s}\right)^{(\rho +d)/2} 
\int_\statset \frac{\psi_0 (y)} {\sqrt{\det M(y)} } \, d\sigma (y).
\end{equation}

We now proceed to complete the proofs by applying \eqref{eq:integrated} in each
case. We give full details only for the case $\b{z^*} = \b{1}$; the general case
again follows from the change of variable $\b{z} \mapsto \b{z}/\b{z^*}$ and we
leave the details to the reader.

\subsection*{Proof of Theorem~\ref{thm:nondeg}}

By the nondegeneracy assumption, $\rank (\cmat) = d+1$, hence 
$\rho := \rank (\overline{\cmat}) = d$ and $\dim (B_y) = 2d$.  
The set $\statset$ has dimension $n - d$; we assume 
this to be strictly positive since the case $n = d$ has been dealt with 
already in Theorem~\ref{thm:compnondeg}. 
At each stationary point $(\b{0} , \bb{\alpha}) \in \statset$ 
we find the Hessian determinant of $f$ restricted to $B_y$ to be the constant
value $(\det \overline{\cmat})^2$.  Now~\eqref{eq:integrated} 
gives
\begin{eqnarray*}
a_\rr & \sim & (2\pi)^{-d} p_0 (s) s^{-d} (2 \pi)^d \phi (\b{1}) 
   \int_\statset \det M(y)^{-1/2}  \, d\sigma (y) \\[1ex]
& \sim & s^{n-d} \frac{ \phi (\b{1}) \sigma (\statset)} {|\det 
   \overline{\cmat}| } \, ,
\end{eqnarray*}
which establishes the formula for the leading term.  That the expansion 
is a polynomial plus an exponentially smaller correction again follows  
from~\cite[Theorem~3.1]{pemantle-finite-dimension}.

\begin{unremark}
The expansion via Leibniz' formula is not generally integrable past
the leading term, and therefore does not give a method of computing
the lower terms of the polynomial $P$. 
\noproof
\end{unremark}

\subsection*{Proof of Theorem~\ref{thm:degenerate}}

Here $\cone_0 = \{ \direc \}$, and so $\mathcal{A}(\direc) = \simp$. Thus
$\dim \mathcal{A}^\perp = 0$, so that for
$y= (\b{0}, \bb{\alpha}) \in \statset$, $M(y) = Q(\b{z^*}, \bb{\alpha})$.
 Use~(\ref{eq:integrated}) once more with $\dim (B_y) = d$ to get 
$$
a_\rr \sim (2\pi)^{-d} p_0 (s) s^{-d/2} (2 \pi)^{d/2} 
 \int_\simp  \frac{\phi (\b{z})}{\sqrt{\det Q(\b{z^*}, \bb{\alpha})} } \, 
   d\mu (\bb{\alpha})
$$
which simplifies to the conclusion of the theorem.
\noproof

\section{Further discussion}
 \tlabel{sec:future}

{\bf Asymptotics in the gaps.} The problem of determining asymptotics
when $\b{\w{r}} / s$ converges to the boundary of $\cone_0$ is dual to
the problem raised in~\cite{pemantle-wilson-multivariate1} of letting
$\b{\w{r}} / s$ converge to $\partial \cone_0$ from the outside. 
Solutions to both of these problems are necessary before we understand
asymptotics ``in the gaps'', that is, in any region asymptotic to and
containing a direction in the boundary of $\cone_0$.  For example, in
the case Example~\ref{eg:tangent2D} with $\direc (1,1)$, what are the
asymptotics for $a_{r , r + \sqrt{r}}$ as $r \to \infty$? 

In addition to the obvious discontinuity in our asymptotic expansions above at the boundary of $\cone_0$ caused by the change from multiple to smooth point regime, in the nondegenerate case other discontinuities arise. The leading term may change degree around the boundary of $\cone_0$. This reflects the behavior of the corresponding Fourier-Laplace integrals, as the affine subspace $\mathcal{A}$ of positive dimension moves from intersecting the interior of $\simp$ to only intersecting  a lower-dimensional face. Subtler effects occur as $\mathcal{A}$ passes a corner in the interior of $\cone$, when the leading term is continuous but lower terms may not be. 
 
A proper analysis of all these issues in the Fourier-Laplace integral 
framework requires the consideration of integrals whose phase 
can vary with $\lambda$.  This is the subject of work in
progress by M. Lladser \cite{lladser-thesis}.


{\bf Effective computation.} The greatest obstacle to making all these computations completely
effective lies in the location of the minimal point $\b{z^*}$ given 
$\rr$.  Assuming the existence of a $\b{z^*} (\rr)$ with 
$\rr \in \cone (\b{z^*})$, how may we compute $\b{z^*} (\rr)$ and
test whether it is a minimal point?  Since the moduli of
the coordinates of $\b{z^*}$ are involved in the definition
of minimality, this is a problem in real rather than complex
computational geometry and does not appear easy.  For example, 
how easily can one prove that the double point in 
Example~\ref{eg:figure 8} is minimal?

{\bf Minimal points.} Many of our theorems rule out analysis of a minimal point
$\b{z^*}$ if one of its coordinates $(z^*)_j$ is zero.  The directions
in $\cone (\b{z^*})$ will always have $r_j = 0$, in which case the
analysis of coefficient asymptotics reduces to a case with 
one fewer variable.  Thus it appears no generality is lost.  
If we are, however, able to solve the previous problem, wherein
$\b{r} / s$ converges to $\partial \cone_0$, then we may choose
to let $\b{r} / s$ converge to something with a zero component.
The problem of asymptotics when some $r_j = o(r_k)$ now makes sense 
and is not reducible to a previous case.  Presumably these asymptotics
are governed by the minimal point $\b{z^*}$ still, but it must be
sorted out which of our results persist when $(z^*)_j = 0$.  Certainly
the geometry near $\b{z^*}$ has more possibilities, since it is easier
to be a minimal point (it is easier to maintain $|z_j| \geq |(z^*)_j|$
for $\b{z}$ near $\b{z^*}$ when $(z^*)_j = 0$).  

In this article and its predecessor we have not treated 
toral minimal points (minimal points $\b{z^*}$ that are not 
strictly minimal, and not isolated in $\torus(\b{z^*})$). 
Nor have we analysed strictly minimal points that are cusps. 
The extension of our results here to those cases is in fact 
rather routine, and will be carried out in future work.

{\bf Local geometry of critical points.} Many of our results fail when the Hessian of $f$  unexpectedly 
vanishes.  Surprisingly, we do not know whether this ever happens 
in cases of interest.  Specifically, we do not know of a generating
function with nonnegative coefficients, for which the Hessian
of $f$ on $B_y$ (notation of Section~\ref{sec:proofs}) ever 
vanishes.  It goes without saying, therefore, that we do not know
whether cases ever arise of degeneracies of mixed orders, such as
may arise in the case of the remark following Example~\ref{eg:tangent2D}.

{\bf Homological methods.} The methods of this paper may be described 
as reasonably elementary though somewhat involved.  The necessary complex  contour and stationary phase integrals are well known, except perhaps for the "resolving lemma", Lemma 4.2.  While writing down these results,  we found a more sophisticated approach involving some algebraic topology (Stratified Morse Theory) and singularity theory.  This approach delivers similar results in a more general framework.  Among other things, it clarifies and generalizes formulae such as the expression~(\ref{eq:P}) derived by hand above.  Still, the present approach has an advantage other than its chronological priority and its independence from theory not well known by practitioners of analytic combinatorics.  Homological methods, as far as we can tell, are not capable of determining asymptotics in "boundary" directions, such as are given in part (ii) of Theorem 3.1. It seems that one needs to delve into the fine details of the integral, as we do in the present work, in order to handle these cases.  Much of the interest in asymptotic enumeration centers on phase boundaries, so the weight of these cases is far from negligible.  We hope some day to integrate the two approaches, but at present each does something the other cannot.

\subsection*{Acknowledgement}
We thank Andreas Seeger for giving a reference to an easier proof 
of Lemma~\ref{lem:simplex}.

\bibliographystyle{alpha}
\bibliography{mvGF}
\end{document}